\documentclass[11pt]{amsart}
\usepackage{verbatim}
\usepackage{amssymb}

\theoremstyle{plain}
\newtheorem{lem}{Lemma}
\newtheorem{thm}[lem]{Theorem}
\newtheorem{prop}[lem]{Proposition}
\newtheorem{cor}[lem]{Corollary}
\newtheorem{example}[lem]{Example}

\begin{document}

\title{$q$ and $q,t$-Analogs of Non-commutative 
Symmetric Functions}
\author{Nantel Bergeron and Mike Zabrocki}
\email{bergeron@mathstat.yorku.ca, 
zabrocki@mathstat.yorku.ca}
\address{Mathematics and Statistics, York University, Toronto, ON
M3J 1P3, Canada}
\maketitle

\font\Ch=msbm10  
\def\scr{\bf}
\def\e{\hbox{\scr e}}
\def\h{\hbox{\scr h}}
\def\s{\hbox{\scr s}}
\def\A{\hbox{\scr A}}
\def\B{\hbox{\scr B}}
\def\H{\hbox{\scr H}}
\def\W{\hbox{\scr W}}
\def\Hw{\widetilde{\H}}
\def\Qs{{F}}
\def\Qm{{M}}
\def\QPq{{P}}
\def\QGq{{G}}
\def\Qsym{{Qsym}}
\def \Q {\hbox{\Ch Q}}
\def \Z {\hbox{\Ch Z}}
\def\NCL{{NC\Lambda}}
\def\la{{\lambda}}
\def\intersect{{\cap}}
\def\AAq{{\widetilde{A_\alpha}}^q}
\def\Axq #1 {{\widetilde{A_{#1}}}^q}
\newcommand{\hatq}[1]{{\widetilde{#1}}^q}
\newcommand{\ob}[1]{{\overline{#1}}}
\def\finereq{{\leq}}
\def\coarsereq{{\geq}}
\def\finer{{<}}
\def\coarser{{>}}
\def\Prop{{\noindent {\bf Proposition}}}
\def\La{{\Lambda}}
\def\unit{{\eta}}
\def\counit{{\varepsilon}}
\def\ncmu{{\mu}}
\def\ncDelta{{\Delta}}
\def\ncunit{{\eta}}
\def\nccounit{{\varepsilon}}
\def\ncantipode{{S}}
\newcommand{\invb}[1]{{\overleftarrow{#1}}}
\newcommand{\invc}[1]{{#1^c}}
\newcommand{\inva}[1]{{#1'}}
\newcommand{\mC}[2]{{\mathcal{C}}_{#2}^#1}
\def\oma{{\inva{\omega}}}
\def\omb{{\invb{\omega}}}
\def\omc{{\invc{\omega}}}
\def\coeff{{\big|}}
\def\tocom{{\chi}}
\def\append{{\cdot}}
\def\attach{{|}}
\def\antipode{{S}}

\edef\savecatcodeat{\the\catcode`@}
\catcode`\@=11

\def\tb@ifSpecChars#1#2{#1}
\def\tb@ifNoSpecChars#1#2{#2}

\def\tableau{%
  \bgroup
  \@ifstar{\let\Tif\tb@ifNoSpecChars\tb@tableauB}
          {\let\Tif\tb@ifSpecChars\tb@tableauB}}

\def\tb@tableauB{
  \@ifnextchar[{\tb@tableauC}{\tb@tableauC[]}}

\def\tb@tableauC[#1]{\hbox\bgroup%
    \let\\=\cr
    \def\bl{\global\let\tbcellF\tb@cellNF}%
    \def\tf{\global\let\tbcellF\tb@cellH}
%
    \dimen2=\ht\strutbox \advance\dimen2 by\dp\strutbox%
    \ifx\baselinestretch\undefined\relax%
    \else%
       \dimen0=100sp \dimen0=\baselinestretch\dimen0%
       \dimen2=100\dimen2 \divide\dimen2 by\dimen0%
    \fi%
    \let\tpos\tb@vcenter
    \tb@initYoung
    \tb@options#1\eoo
    \let\arrow\tb@arrow%
    \dimen0=\Tscale\dimen2%
    \dimen1=\dimen0 \advance\dimen1 by \tb@fframe%
    \lineskip=0pt\baselineskip=0pt
%
    \def\tb@nothing{}%
    \def\endcellno{$\rss\egroup\bss\egroup}
    \def\endcell{\endcellno\kern-\dimen0}
    \def\begincell{\vbox to\dimen0\bgroup\vss\hbox to\dimen0\bgroup\hss$}%
    \let\overlay\tb@overlay%
    \let\fl\tb@fl%
    \let\lss\hss\let\rss\hss\let\tss\vss\let\bss\vss
    \def\mkcell##1{
        \let\tbcellF\tb@cellD
        \def\tb@cellarg{##1}
        \ifx\tb@cellarg\tb@nothing\let\tb@cellarg\tb@cellE\fi%
        \begincell\tb@cellarg\endcellno
        \tbcellF}
    \let\savecellF\tbcellF
     \Tif{\catcode`,=4\catcode`|=\active}{}\tb@tableauD}%

\let\tb@savetableauD\tableauD
{
    \catcode`|=\active \catcode`*=\active \catcode`~=\active%
    \catcode`@=\active
\gdef\tableauD#1{%
  \Tif{
    \mathcode`|="8000 \mathcode`*="8000%
    \mathcode`~="8000 \mathcode`@="8000%
    \def@{\bullet}%
    \let|\cr
    \let*\tf
    \let~\sk
  }{}%
  \tpos{\tabskip=0pt\halign{&\mkcell{##}\cr#1\crcr}}%
  \global\let\tbcellF\savecellF
  \egroup
  \egroup}
}
\let\tb@tableauD\tableauD
\let\tableauD\tb@savetableauD
\let\tb@savetableauD\undefined


\def\tb@options#1{\ifx#1\eoo\relax\else\tb@option#1\expandafter\tb@options\fi}

\def\tb@option#1{%
  \if#1t\let\tpos\tb@vtop\fi
  \if#1c\let\tpos\tb@vcenter\fi
  \if#1b\let\tpos\vbox\fi
  \if#1F\tb@initFerrers\fi
  \if#1Y\tb@initYoung\fi
  \if#1s\tb@initSmall\fi
  \if#1m\tb@initMedium\fi
  \if#1l\tb@initLarge\fi
  \if#1p\tb@initPartition\fi
  \if#1a\tb@initArrow\fi
}

\def\tb@vcenter#1{\ifmmode\vcenter{#1}\else$\vcenter{#1}$\fi}

\def\tb@vtop#1{\hbox{\raise\ht\strutbox\hbox{\lower\dimen0\vtop{#1}}}}

\def\tb@initPartition{\def\Tscale{.3}}
\def\tb@initSmall{\def\Tscale{1}}
\def\tb@initMedium{\def\Tscale{2}}
\def\tb@initLarge{\def\Tscale{3}}

\def\tb@initArrow{\dimen2=1.25em}

\def\tb@initYoung{%
  \def\tb@cellE{}
  \let\tb@cellD\tb@cellN
  \def\sk{\global\let\tbcellF\tb@cellNF}}
\def\tb@initFerrers{%
  \def\tb@cellE{\bullet}
  \let\tb@cellD\tb@cellNF
  \def\sk{\bullet}}

\tb@initMedium

\def\tb@sframe#1{%
  \vbox to0pt{
    \vss
    \hbox to0pt{%
      \hss
      \vbox to\dimen1{
        \hrule depth #1 height0pt
        \vss
        \hbox to\dimen1{
          \vrule width #1 height\dimen1
          \hss
          \vrule width #1
          }%
        \vss
        \hrule height #1 depth 0in
        }%
      \kern-\tb@hframe
      }%
    \kern-\tb@hframe}}

\def\tb@hframe{.2pt}\def\tb@fframe{.4pt}\def\tb@bframe{2pt}
\def\tb@cellH{\tb@sframe{\tb@bframe}}       
\def\tb@cellNF{}                            
\def\tb@cellN{\tb@sframe{\tb@fframe}}       
\let\tbcellF\tb@cellN                       

\def\tb@rpad{1pt}
\def\tb@lpad{1pt}
\def\tb@tpad{1.8pt}
\def\tb@bpad{1.8pt}

\def\tb@overlay{\endcell\@ifnextchar[{\tb@overlaya}{\begincell}}
\def\tb@overlaya[#1]{\vbox to\dimen0\bgroup%
  \tb@overlayoptions#1\eoo%
  \tss\hbox to\dimen0\bgroup\lss$}
\def\tb@overlayoptions#1{\ifx#1\eoo\relax\else\tb@overlayoption#1\expandafter\tb@overlayoptions\fi}

\def\tb@overlayoption#1{
  \if#1t\def\tss{\vskip\tb@tpad}\let\bss\vss\fi
  \if#1c\let\tss\vss\let\bss\vss\fi
  \if#1b\def\bss{\vskip\tb@bpad}\let\tss\vss\fi
  \if#1l\def\lss{\hskip\tb@lpad}\let\rss\hss\fi
  \if#1m\let\lss\hss\let\rss\hss\fi
  \if#1r\def\rss{\hskip\tb@rpad}\let\lss\hss\fi
}

\def\tb@fl{\endcell\begincell\vrule depth 0pt width \dimen0 height \dimen0 \endcell\begincell}



\@ifundefined{diagram}{}{
\def\tb@arrowpad{.5}

\newoptcommand{\tb@arrow}{\@ne}[2]{%
  \endcell
   \begingroup%
   \let\dg@getnodesize\tb@getnodesize
   \dg@USERSIZE=#1\relax%
   \ifnum\dg@USERSIZE<\@ne \dg@USERSIZE=\@ne \fi%
   \dg@parse{#2}%
   \dg@label{\tb@draw{#1}{#2}}}

\def\tb@getnodesize#1#2#3#4#5{\dimen3=\tb@arrowpad\dimen2 #4=\dimen3 #5=\dimen3\relax}
\def\tb@getnodesize#1#2#3#4#5{\ifnum#2=0\ifnum#3=0\tb@getnodesizetail{#4}{#5}\else\tb@getnodesizehead{#4}{#5}\fi\else\tb@getnodesizehead{#4}{#5}\fi}
\def\tb@getnodesizetail#1#2{\dimen3=.5\dimen2 #1=\dimen3 #2=\dimen3}
\def\tb@getnodesizehead#1#2{\dimen3=.5\dimen2 #1=\dimen3 #2=\dimen3}

\def\tb@draw#1#2#3#4{%
        \dg@X=0\dg@Y=0\dg@XGRID=1\dg@YGRID=1\unitlength=.001\dimen0%
        \dg@LBLOFF=\dgLABELOFFSET \divide\dg@LBLOFF\unitlength%
        \dg@drawcalc
        \begincell
        \let\lams@arrow\tb@lams@arrow
        \begin{picture}(0,0)\begingroup\dg@draw{#1}{#2}{#3}{#4}\end{picture}%
        \endcell
        \endgroup
        \begincell}
}

%
%
%
\def\tb@lams@arrow#1#2{%
 \lams@firstx\z@\lams@firsty\z@
 \lams@lastx#1\relax\lams@lasty#2\relax
 \lams@center\z@
 %
 \N@false\E@false\H@false\V@false
 \ifdim\lams@lastx>\z@\E@true\fi
 \ifdim\lams@lastx=\z@\V@true\fi
 \ifdim\lams@lasty>\z@\N@true\fi
 \ifdim\lams@lasty=\z@\H@true\fi
 \NESW@false
 \ifN@\ifE@\NESW@true\fi\else\ifE@\else\NESW@true\fi\fi
 %
 \ifH@\else\ifV@\else
  \lams@slope
  \ifnum\lams@tani>\lams@tanii
   \lams@ht\ten@\p@\lams@wd\ten@\p@
   \multiply\lams@wd\lams@tanii\divide\lams@wd\lams@tani
  \else
   \lams@wd\ten@\p@\lams@ht\ten@\p@
   \divide\lams@ht\lams@tanii\multiply\lams@ht\lams@tani
  \fi
 \fi\fi
 %
 \ifH@  \lams@harrow
 \else\ifV@ \lams@varrow
 \else \lams@darrow
 \fi\fi
}

\catcode`\@=\savecatcodeat
\let\savecatcodeat\undefined


\begin{abstract}
We introduce two families of non-commutative
symmetric functions that have analogous properties
to the Hall-Littlewood and Macdonald symmetric
functions.
\end{abstract}

\section{Introduction}
 It was noticed in
\cite{Z} that many $q$ and $q,t$ analogs that are commonly
studied in the space of symmetric functions
arise from an unusual $q$-twisting of the
symmetric function found by setting $q=0$.  In particular
we see that the Hall-Littlewood and Macdonald \cite{M}
symmetric functions arise from this construction by taking a
$q$-analog of operators that add a row on the Schur basis
or a column on the Hall-Littlewood symmetric functions.

The $q$-twisting that was given in that article may
easily be expressed in terms of the notation of Hopf
algebras (since the symmetric functions form a commutative
and co-commutative Hopf algebra structure). 
Within this context it can be seen that this
may be generalized to any graded Hopf
algebra.  It is not clear when or if this $q$-analog
will be interesting in another setting, but it creates a 
context for looking
for `natural' examples of $q$-analogs within other spaces.

The non-commutative symmetric functions \cite{GKLLRT} form 
such a graded
Hopf algebra and are an obvious place to begin searching
for such examples.  The natural analog of the Schur functions
within this space are the ribbon Schur functions.  

In a manner analogous to that in the symmetric functions,
it is possible to define operators that add a row to the
ribbon Schur functions.  The $q$-analog of these operators
are a natural analog of the operators
that add a row to the Hall-Littlewood symmetric functions.
Remarkably, we see that this action gives rise to
a family of non-commutative symmetric functions
that when expanded in terms of ribbon Schur functions
have coefficients that are a power of $q$.  

These NC-symmetric functions have properties that suggest
they are an analog of the
Hall-Littlewood symmetric functions,
but they are not equivalent to the non-commutative Hall-Littlewood
symmetric functions of Hivert \cite{Hi}.
They do however share some of the same properties of
Hivert's NC-Hall-Littlewoods including a factorization
property by setting $q$ equal to a root of unity
(Proposition \ref{factor}).  It
should be remarked that the differences
between the factorization properties of these and the
Hivert functions
suggest that the Hivert functions are
associated to the length statistic 
on compositions in the
same way that the Hall-Littlewood functions presented
here are associated to the size statistic
on compositions.

This family of NC-symmetric functions 
$\{\H_\alpha^q\}_\alpha$ have the following distinguishing
properties.

\vskip .1in
\noindent
1. They are triangularly related to the ribbon Schur basis.
Namely, we have 
\begin{equation*}
\H_\alpha^q = \s_\alpha + \sum_{\beta>\alpha}
c^q_{\alpha\beta} \s_\beta.
\end{equation*}

\vskip .1in
\noindent
2. The coefficient of a single ribbon function in 
$\H_\alpha^q$ is either $0$ or a power of $q$.
The coefficient of the ribbon function indexed by
a single part in $\H_\alpha^q$ is $q^{n(\alpha)}$ where
$n(\alpha) = \sum_i i$ where the sum is over all $i$ 
which are descents of $\alpha$.

\vskip .1in
\noindent
3. When $q=1$, $\H_\alpha^q$ becomes $\h_\alpha$,
the non-commutative analogs of the homogeneous symmetric
functions.  When $q=0$,
$\H_\alpha^q$ specializes to the ribbon Schur function 
$\s_\alpha$.  When $q$ is a root of unity then $\H_\alpha^q$
specializes to a product of non-commutative symmetric
functions.

Considering the morphism $\chi$ that sends the non-commutative
symmetric function $\h_\alpha$ to the commutative version
$h_\alpha$, the image of the functions $\H_\alpha^q$ are
the commutative Hall-Littlewood functions whenever
the composition $\alpha$ represents a
partition (i.e. when $\alpha$ is a hook).  That is,

\vskip .1in
\noindent
4. (Proposition \ref{Hqcomm}) $\chi( \H_{(1^a,b)}^q) = H_{(b,1^a)}^q$ where $H_\la^q =
\sum_\mu K_{\mu\la}(q) s_\mu$.

We introduce an inner product on the space of
non-commutative symmetric functions by setting the ribbon
Schur functions as `semi-self' dual in the following manner
\begin{equation*}
\left< \s_\alpha, \s_\beta \right> = 
(-1)^{|\alpha|+\ell(\alpha)} \delta_{\alpha\invc{\beta}}.
\end{equation*}

This inner product does
not seem to appear elsewhere in the literature, but does
share similar properties with the inner product of the
symmetric functions and can be a useful tool for 
calculation within this space.  The surprising property
that we observe is that the non-commutative analogs of
the elementary, homogeneous and Hall-Littlewood bases also
share this `semi-self' duality property.  That is, we have
in addition to the properties mentioned above,

\vskip .1in
\noindent
5. (Proposition \ref{Hscalar})
\begin{equation*}
\left< \H_\alpha^q, \H_\beta^q \right> = 
(-1)^{|\alpha|+\ell(\alpha)} \delta_{\alpha\invc{\beta}}.
\end{equation*}

Most of these properties are analogous to ones
that exist for the 
non-commutative Hall-Littlewood analogues of Hivert
\cite{Hi}, however this last property is
not shared by Hivert's noncommutative Hall-Littlewood
functions.

Next, we consider a $q,t$-analog of the non-commutative
symmetric functions where we look for properties that
are analogous to the Macdonald symmetric functions.  
Of course, many $q,t$-analogues are possible and we
consider one that that has properties which generalize
those for our version of the non-commutative 
Hall-Littlewood and seem to be analogous to the
Macdonald symmetric functions.  The family 
$\{\H_\alpha^{qt}\}_\alpha$ has the following
important properties:

\vskip .1in
\noindent
1. There is a triangular relation between the family
$\{ \H_\alpha^t\}_\alpha$ and $\{ \H_\alpha^{qt}\}_\alpha$.
\begin{equation*}
\H_\alpha^{qt} = 
\sum_{\beta \leq \alpha}
c_{\alpha\beta}^{qt} \H_\beta^{t}.
\end{equation*}

\vskip .1in
\noindent
2. The coefficient of a single ribbon function in 
$\H_\alpha^{qt}$ is of the form $q^a t^b$ (with
$a,b\geq 0)$.
The coefficient of a ribbon Schur function
indexed by a single part in
$\H_\alpha^{qt}$ is $t^{n(\alpha)} $, the coefficient of
a ribbon Schur function indexed by a composition
of $1$s is $q^{n(\inva{\alpha})}$.

\vskip .1in
\noindent
3. We have the specialization
$\H_\alpha^{0t} = \H_\alpha^t$, and
$\H_\alpha^{1t}$ is a product of
non-commutative symmetric functions
(Proposition \ref{H1tprod}).

\vskip .1in
\noindent
4. The $\H_\alpha^{qt}$ satisfy the following two relations
\begin{equation*}
\H_\alpha^{tq} = \oma \H_{\inva{\alpha}}^{qt},
\end{equation*}
\begin{equation*}
q^{n(\inva{\alpha})} t^{n(\alpha)} 
\H_\alpha^{\frac{1}{q}\frac{1}{t}} =
\omc \H_{\alpha}^{qt}.
\end{equation*}

\vskip .1in
\noindent
5. (Proposition \ref{Hqtcomm}) 
$\chi( \H_{(1^a,b)}^{qt}) = H_{(b,1^a)}^{qt},$ 
where $H_\la^{qt} = \sum_\mu K_{\mu\la}(q,t) s_\mu$.

\vskip .1in
\noindent
6. (Proposition \ref{Hqtscalar})
\begin{equation*}
\left< \H_{\alpha}^{qt}, \H_{\beta}^{qt} \right>
= (-1)^{|\alpha| + \ell(\alpha)}
\delta_{\alpha\invc{\beta}}
\prod_{i=1}^{n-1} (1-q^i t^{n-i}) .
\end{equation*}

The most remarkable property that arises from these functions
is the existence of an operator $\nabla$ that has the
NC-Macdonald functions as eigenfunctions.  That is, if we
set
\begin{equation*}
\nabla \Hw_\alpha^{qt} = q^{n(\inva{\alpha})} t^{n(\alpha)}
\Hw_\alpha^{qt},
\end{equation*}
where $\Hw_\alpha^{qt} = t^{n(\alpha)} \H_\alpha^{q\frac{1}{t}}$,
then this operator has an elegant action on
the ribbon Schur functions and 
shares many of the same properties that exist in the 
commutative case \cite{BGHT}.  Unlike in the commutative case
however, formulas for this operator are immediately
solvable.  In a beautiful analogy, where
the commutative version of the operator
$\nabla$ produced a $q,t$ grading of the space of parking
functions, the non-commutative $\nabla$ produces a grading
of the space of preferential arrangements.

In searching for an interesting non-commutative analog
of the Macdonald symmetric functions, we considered many
possibilities (including the analog considered in 
\cite{HLT}).  None of the functions except for the one
we discuss here had an equivalent $\nabla$ operator and 
it was this property that indicated to us that these
functions are indeed remarkable.

It is not known yet if this family has a
representation theoretical model analogous to
the $n!$-conjecture or the diagonal harmonics
that motivate the existence of these functions.
We do see however that the
non-commutative versions of these functions and operators
share many of the same properties with the commutative case.
Independent of their own interest,
it is at least hopeful that they will give some insight 
into why the some of the conjectures for the commutative
case are true.

\section{Notation for compositions, partitions, Hopf algebras,
symmetric, NC-symmetric, and Quasi-symmetric functions}
\subsection{Compositions} \label{comps}

We will say that $\alpha$ is a composition of $n$ and write
$\alpha \models n$ if
$\alpha$ is a sequence of positive integers such that
$\alpha_1 + \alpha_2 +
\cdots + \alpha_{\ell(\alpha)} = |\alpha| = n$. The length of
the sequence is denoted by the symbol
$\ell(\alpha)$.

For any two compositions $\alpha$ and $\beta$, define
the {\it concatenate} and the {\it attach} of $\alpha$ and $\beta$
to be the 
compositions (respectively) 
\begin{equation}
\alpha \append \beta = (\alpha_1, \alpha_2,
\ldots, \alpha_{\ell(\alpha)}, \beta_1, \beta_2,
\ldots, \beta_{\ell(\beta)})
\end{equation}
and
\begin{equation}
\alpha \attach \beta = (\alpha_1, \alpha_2,
\ldots, \alpha_{\ell(\alpha)} + \beta_1, \beta_2,
\ldots, \beta_{\ell(\beta)}).
\end{equation}

For a composition $\alpha$ of $n$ define the descent set of
$\alpha$ to be the subset of $\{ 1,2,3, \ldots, n-1\}$ as
$D(\alpha) = \{\alpha_1, \alpha_1+\alpha_2, \ldots,
\alpha_1 +\alpha_2+ \cdots + \alpha_{\ell(\alpha)-1}\}$.
The size of the descent set of $\alpha$ is one less than
it's length and it is easily seen that the compositions
of $n$
are in one-to-one correspondence with the subsets of
$\{1, 2, \ldots, n-1\}$.

There  is a natural partial order on the compositions
of $n$.  Say that a composition $\alpha$ is finer than
a composition $\beta$ (or $\beta$ is coarser than $\alpha$)
and write $\alpha \finereq \beta$
if there exists compositions $\gamma^{(1)}, \gamma^{(2)},
\ldots, \gamma^{(k)}$ such that $\alpha = 
\gamma^{(1)} \append \gamma^{(2)} \append
\ldots \append \gamma^{(k)}$ and $\beta =
\gamma^{(1)} \attach \gamma^{(2)} \attach
\ldots \attach \gamma^{(k)}$.  Alternatively, in terms of
descent sets we say that $\alpha \finereq \beta$ if and
only if $D(\beta) \subseteq D(\alpha)$.

There are three standard involutions on the set
of compositions.  The first involution reverses the
order of the sequence.  
We set $\invb{\alpha} = (\alpha_{\ell(\alpha)},
\alpha_{\ell(\alpha)-1}, \ldots,\\ \alpha_{1})$.  If the
descent set of $\alpha$ is $D(\alpha) = \{ i_1, i_2, 
\ldots, i_k \}$ then $D(\invb{ \alpha}) = 
\{ |\alpha| - i_1, |\alpha|-i_2,
 \ldots, |\alpha|-i_k \}$.
 
The second involution corresponds to taking the complement
of the descent set.  Define
$\invc{\alpha}$
to be the composition
with $D({\invc{\alpha}}) = \{ 1, 2, \ldots, |\alpha| - 1 \} -
D(\alpha)$.  Notice that if $\alpha$ is a composition
of length $k$ then 
$\invc{\alpha}$ 
is a composition of length
$n+1-k$.

Finally, the third involution is the composition of the
previously two defined.  Let $\inva{\alpha} = 
\invb{\invc{\alpha}}=
\invc{\invb{\alpha}}$. 
The composition 
$\inva{\alpha}$
also has length $n+1-k$.

\begin{figure}
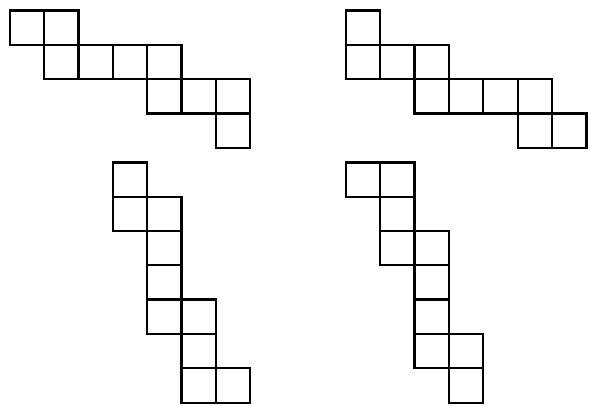

$$
\tableau[s]{,|\bl,,,,|\bl,\bl,\bl,\bl,,,|\bl,\bl,\bl,\bl,\bl,\bl,|}\hskip .5in
\tableau[s]{|,,|\bl,\bl,,,,|\bl,\bl,\bl,\bl,\bl,,}$$
$$\tableau[s]{|,|\bl,|\bl,|\bl,,|\bl,\bl,|\bl,\bl,,}\hskip .5in
\tableau[s]{,|\bl,|\bl,,|\bl,\bl,|\bl,\bl,|\bl,\bl,,|\bl,\bl,\bl,|}$$
\caption{Images of the three involutions.  
If $\alpha = (2,4,3,1)$, 
then ${\overline{\alpha}} = (1,3,4,2)$, 
$\invc{\alpha} = (1,2,1,1,2,1,2)$
and  $\inva{\alpha} = (2,1,2,1,1,2,1)$ .}
\end{figure}

For some formulas we will need a total order on the 
compositions of size $n$.   We
set $\phi(\alpha) = \sum_{i \in D(\alpha)} 2^{i-1}$.
This map associates each composition with an integer
between $0$ and $2^{n-1} - 1$.  This map induces a total
order from the integers on this set which is a refinement
of the partial order defined above.

For the composition $\alpha$, there is a standard
statistic we will use frequently given by $n(\alpha)
= \sum_{i \in D(\alpha)} i$.

\subsection{Partitions}

A partition $\la$ of $n$ is a composition of $n$ with the
property that $\la_1 \geq \la_2 \geq \cdots \geq
\la_{\ell(\la)}$.  We will indicate that $\la$ is a partition
by the notation that $\la \vdash n$.  We say that a
partition $\mu$ is contained in a partition $\la$ and write
that $\mu \subseteq \la$ if 
$\ell(\mu) \leq \ell(\la)$ and
$\mu_i \leq \la_i$ for all $1 \leq i \leq \ell(\mu)$.
Define then a skew partition to be represented by
$\la \slash \mu$ where $\la$ and $\mu$ are partitions
such that $\mu \subseteq \la$.

There is a partial order on the set of partitions.  We
will say that the partition $\la \geq \mu$ if $\la_1 +
\la_2 +\cdots+ \la_i \geq \mu_1 + \mu_2 + \cdots +\mu_i$
for all $i \geq 1$.

Partitions will sometimes be represented by Ferrers diagrams,
a graphical representation of a partition formed by placing
rows of square cells aligned on the left hand edge with $\la_i$
cells in the $i^{th}$ row.  We will use the cartesian convention
where we place the $1^{st}$ row of cells on the bottom
of the diagram (the matrix convention is to place the $1^{st}$
row of cells on the top).  
A skew Ferrers diagram for a skew partition
$\la\slash\mu$
is Ferrers diagram for the partition $\la$ where the cells
that correspond to the partition $\mu$ are not drawn.

Define the conjugate partition to $\la$ to be the partition
$\la'$ such that $\la'_i$ is the number of parts of $\la$
that have size greater than or equal to $i$.
This corresponds to the
partition formed by flipping $\la$ across the line
$x=y$.

Any composition of $n$ may be associated with a `ribbon,'
a skew partition with no $2 \times 2$ sub-diagrams.  This
ribbon is usually represented by a skew-Ferrers diagram. The
composition $\alpha \vdash n$ is mapped to the skew partition
$(\alpha_1 + \alpha_2 + \cdots + \alpha_{\ell(\alpha)} -
\ell(\alpha) + 1,
\alpha_1 + \alpha_2 + \cdots + \alpha_{\ell(\alpha)-1} -
\ell(\alpha) + 2, \ldots, \alpha_1) \slash
(\alpha_1 + \alpha_2 + \cdots + \alpha_{\ell(\alpha)-1} -
\ell(\alpha)+1,
\alpha_1 + \alpha_2 + \cdots + \alpha_{\ell(\alpha)-2} -
\ell(\alpha)+2, \ldots, \alpha_1+\alpha_2-2,\alpha_1-1)$.
For the example in Figure 1, $\alpha = (2,4,3,1)$ is
associated to the skew partition $(7,7,5,2)/(6,4,1)$.

We will label the cells of the $x,y$-coordinate
lattice and
say that a point $(i,j)$ is in the diagram of a partition
$\mu$ if $1 \leq i \leq \mu_j$.  The arm of the cell
$s=(i,j)$ in a partition $\mu$ is denoted by $a_\mu(s) :=
\mu_j - i$.  The leg will be denoted by the value $l_\mu(s)
:= a_{\mu'}(j,i)$.

\subsection{Hopf algebras}
For general facts about Hopf algebras, we refer the
reader to \cite{A} or \cite{S}.

Let $R$ be a commutative ring and $H$ an $R$ module.
We say that $H$ is an algebra if there are maps  
$\mu: H \otimes H \rightarrow H$ (multiplication)
and $\unit: R \rightarrow H$ (unit) that satisfy the 
following two conditions:

\noindent
1) $\mu \circ (id \otimes \mu) = \mu \circ (\mu \otimes
id)$ 

\noindent
2) $\mu \circ ( \unit \otimes id)  = id
= \mu \circ ( id \otimes \unit)$

We say that $H$ has a co-algebra structure
if there are maps 
$\Delta: H \rightarrow H \otimes H$ (comultiplication)
and $\counit: H \rightarrow R$ (counit)
that satisfy the
following two conditions:

\noindent
1) $ (\Delta \otimes id) \circ \Delta = (id \otimes \Delta)
\circ \Delta$

\noindent
2) $(\counit \otimes id ) \circ \Delta  = id
= (id \otimes \counit) \circ \Delta$

If $H$ has at the same time an algebra and a co-algebra 
structure $(H,\mu,\unit)$ and $(H,\Delta,\counit)$ and 
$\Delta$
is a homomorphism of algebras then
we call $H$ together with these corresponding operations,
$(H,\mu,\unit,\Delta,\counit)$, a bialgebra.   
$H$ is called a Hopf algebra if
it is a bialgebra with
a map 
$\antipode: H \rightarrow H$ that satisfies the following
identity:

\noindent
1) $\mu \circ \antipode \otimes id\circ \Delta=
\mu \circ id \otimes \antipode \circ \Delta=
\unit \circ\counit$

If we define the map $\tau: H \otimes H \rightarrow H \otimes
H$ by $\tau(a \otimes b) = b \otimes a$ then we say that
an algebra $H$ is commutative if $\mu \circ \tau = \mu$
and we say that a co-algebra $H$ is co-commutative
if $\tau \circ \Delta = \Delta$.  It may be shown that for
any Hopf algebra $\Delta \circ \antipode = \antipode \otimes
\antipode \circ \tau \circ \Delta$.

If the Hopf algebra is either commutative or co-commutative,
then it follows that $\antipode$ is an involution.

An important operation that arises in this setting
is the convolution of two operators $f,g \in Hom(H,H)$.
We set
$f \ast g = \mu \circ
f \otimes g \circ \Delta$.  Convolution is an associative
binary operation and the element $\unit\counit$
serves as the identity.  That is, we have for
$V \in Hom(H,H)$
\begin{equation}
\unit\counit \ast V = V \ast \unit\counit = V.
\end{equation}
In addition, it follows from the defining property of
the antipode that $id \ast \antipode 
= \antipode \ast id = \unit\counit$.

\subsection{Symmetric functions} We refer the reader to
\cite{M} for basic facts about the symmetric functions.

Consider the space of
symmetric functions as the polynomial ring
$\La = \Q[p_1, p_2, \ldots]$
in the commuting set of variables $\{ p_1, p_2, p_3, 
\ldots\}$.  The $p_i$ are the simple power symmetric
functions and represent the symmetric formal series $p_k =
x_1^k + x_2^k + x_3^k + \cdots$ (although in this context
we need not consider the variables 
$\{x_1, x_2, x_3, \ldots\}$).  Define the degree of the
power symmetric function $p_k$ within this space to have
degree $k$.  Let $\La^n$ represent the subspace of
polynomials of degree $n$.  Since the $p_k$ commute, we
see that $\La^n$ is spanned by
the set of monomials $p_\la :=
p_{\la_1} p_{\la_2} \cdots p_{\la_\ell(\la)}$
where $\la$ is a partition of $n$.

Let $n_i(\la)$ represent the number of parts of size
$i$ in the partition $\la$, then define $z_\la =
\prod_{i \geq 1} i^{n_i(\la)} n_i(\la)!$.
The simple elementary symmetric functions  are
defined to be $e_k = \sum_{\la \vdash k} (-1)^{k - \ell(\la)}
p_\la / z_\la$ and the simple homogeneous symmetric
functions are defined to be
$h_k = \sum_{\la \vdash k} p_\la / z_\la$.  Set $e_\la =
e_{\la_1} e_{\la_2} \cdots e_{\la_{\ell(\la)}}$ and
$h_\la = h_{\la_1} h_{\la_2} \cdots h_{\la_{\ell(\la)}}$.
The Schur symmetric functions are defined to be
$s_{\la} = det|h_{\la_i + i -j}| = det|e_{\la_i' + i -j}|$.

There is a natural scalar product defined on this space
that is defined on the power symmetric functions by
$\left< p_\la, p_\mu \right> = z_{\la} \delta_{\la\mu}$.  It may
be shown that the Schur functions are self-dual with
respect to the scalar product, that is, $\left< s_\la, 
s_\mu \right> = \delta_{\la\mu}$.

We are interested in finding analogs for the two
families of symmetric functions
$H_\mu^q := \sum_{\la\vdash|\mu|} K_{\la\mu}(q) s_\la$
and $H_\mu^{qt} := \sum_{\la\vdash|\mu|} K_{\la\mu}(q,t) s_\la$
(the Hall-Littlewood and Macdonald symmetric functions).
For a definition of $K_{\la\mu}(q)$ and $K_{\la\mu}(q,t)$
and their associated properties, we refer the interested
reader to \cite{M}.

\subsection{Non-commutative symmetric functions}
For a more detailed reference about the non-commutative
symmetric functions, we refer the reader to \cite{GKLLRT}.

Consider the space of non-commutative symmetric functions
as the polynomial ring $\NCL = \Q< \h_1, \h_2, \ldots >$
in the non commuting set of variables
$\{ \h_1, \h_2, \h_3, \ldots \}$ .  The degree of a
monomial $\h_{i_1} \h_{i_2} \cdots \h_{i_\ell}$ will be
the sum of the indices $i_1 + i_2 + \cdots +i_\ell$. 
The span of the monomials of $\NCL$ of degree $n$ will be
denoted by $\NCL^n$ so that $\NCL = \bigoplus_{n \geq 0}
\NCL^n$ is a graded ring.

We will define
$\e_k = \sum_{\alpha \models k} (-1)^{k-\ell(\alpha) }
\h_\alpha$.  These are the analogs of the
basic homogeneous and elementary symmetric functions. 
For any composition we define $\h_\alpha =
\h_{\alpha_1}
\h_{\alpha_2}
\cdots \h_{\alpha_{\ell(\alpha)}}$ and $\e_\alpha =
\e_{\alpha_1}
\e_{\alpha_2}
\cdots \e_{\alpha_{\ell(\alpha)}}$.  The sets
$\{\h_\alpha\}_{\alpha\models n}$ and
$\{\e_\alpha\}_{\alpha\models n}$ are all bases for the space
of non-commutative symmetric functions of degree $n$.

The
ribbon Schur functions are defined to be 
$\s_\alpha = \sum_{\beta \coarsereq \alpha}
(-1)^{\ell(\alpha)-\ell(\beta)} \h_\beta$.  It is well known that the set
$\{\s_\alpha\}_{\alpha \models n}$ also defines a basis
for $\NCL^n$.  It is normal to define a pairing of the
NC-symmetric functions with the Quasi-symmetric functions.
Instead here we define a scalar product on this space such
that the ribbon Schur functions are self-dual, that is,
\begin{equation}\left< \s_\alpha, \s_{\beta} \right> = 
(-1)^{|\alpha|+\ell(\alpha)} \delta_{\alpha\invc{\beta}}.
\label{NCscalar}
\end{equation}

We remark that this product has the property that
$\left< f, g \right> 
= (-1)^{deg(f)+1} \left< g, f \right>$.  This follows
from the relation  $\ell(\alpha)+
\ell(\alpha^c) = |\alpha|+ 1$.
We will prove in full generality in later sections (although
it is not difficult to show) that we
have the following relations

\begin{prop}
\begin{equation}\label{selfdual}
\left< \s_\alpha, \s_{\beta} \right> = 
\left< \h_\alpha, \h_\beta \right> =
\left< \e_\alpha, \e_\beta \right> =
(-1)^{|\alpha|+\ell(\alpha)} \delta_{\alpha\invc{\beta}}.
\end{equation}
\end{prop}

\begin{cor} \label{scalarexpr}
For any $f \in \NCL^n$, we have
\begin{equation}
f=\sum_{\alpha \models n} (-1)^{\ell(\alpha)+1}
\left< \s_{\invc{\alpha}}, f \right> \s_\alpha
=\sum_{\alpha \models n} (-1)^{|\alpha|+ \ell(\alpha)}
\left< f, \s_{\invc{\alpha}} \right> \s_\alpha.
\end{equation}
\end{cor}

A similar statement can be made for the $\h_\alpha$
and $\e_\alpha$ bases.

There is a co-commutative Hopf algebra structure on the
space of non-commutative symmetric functions.  
Multiplication $\ncmu$,
the unit $\ncunit$ and the counit $\nccounit$
as defined in the usual manner.  The comultiplication
is defined as the algebra homomorphism that sends
$\ncDelta(\h_k) = 
\sum_{i=0}^k \h_i\otimes\h_{k-i}$ and $\ncDelta(\e_k)
= \sum_{i=0}^k  \e_i\otimes\e_{k-i}$. 
The antipode is defined
by $\ncantipode( \s_\alpha) = (-1)^{|\alpha|} \s_{
\inva{\alpha}}$,
$\ncantipode( \h_\alpha ) = 
(-1)^{|\alpha|} \e_{\invb{\alpha}}$ and
$\ncantipode( \e_\alpha ) = 
(-1)^{|\alpha|} \h_{\invb{\alpha}}$.
Elementary properties
of the antipode and the scalar product show that
$\left< \ncantipode(f), \ncantipode(g) \right> = 
 \left< g, f \right>$.

There are three standard involutions that correspond
to those that exist for the compositions.  Set 
$\oma(\s_\alpha) = \s_{\inva{\alpha}}$, $\omb(\s_\alpha)=
\s_{\invb{\alpha}}$ and $\omc(\s_\alpha) =
\s_{\invc{\alpha}}$.  It is important to note
that we have the relations $\oma \omb = \omb \oma = \omc$.
These operations may be expressed on other bases as well
yielding the following expressions:
\begin{align}\nonumber
\oma(\h_\alpha) = \e_{\invb{\alpha}}\hskip .5in &
\oma(\e_\alpha) = \h_{\invb{\alpha}}
\\\omb(\h_\alpha) = \h_{\invb{\alpha}}\hskip .5in &
\omb(\e_\alpha) = \e_{\invb{\alpha}}
\\\omc(\h_\alpha) = \e_{\alpha}\hskip .5in &
\omc(\e_\alpha) = \h_{\alpha} \nonumber
\end{align}
Of course we also see that $\left< \omb{f}, \omb{g}\right>
= \left< f,g \right>$ and $\left<\oma{f}, \oma{g} \right> =
\left<\omc{f}, \omc{g} \right> = \left< g,f
\right>$.

We will sometimes wish to look at the commutative versions
of the non-commutative symmetric functions.  To this end,
we introduce the surjection 
$\tocom: \NCL \rightarrow \La$
which sends $\h_\alpha$ to the symmetric function $h_{\alpha_1} 
h_{\alpha_2} \ldots h_{\alpha_k}$.  

\subsection{The quasi-symmetric functions}
Consider the space of polynomials in the commuting
set of variables $x_1, x_2, x_3, \ldots, x_n$.  The 
quasi-symmetric functions will be denoted by $\Qsym$
which will be the subspace of polynomials spanned by
the functions
\begin{equation}
\Qm_\alpha = \sum_{f} x_{f(1)}^{\alpha_1} x_{f(2)}^{\alpha_2}
\cdots x_{f(\ell(\alpha))}^{\alpha_{\ell(\alpha)}},
\end{equation}
where the sum is over all functions $f: [\ell(\alpha)]
\rightarrow [n]$ such that $f(i) < f(i+1)$.

These functions are the analogs of the monomial symmetric
functions within the space of symmetric functions.
There is a standard pairing between the quasi-symmetric
functions and space of non-commutative symmetric functions.
This pairing is defined by setting non-commutative
homogeneous symmetric functions as dual to the $\Qm_\alpha$
basis, that is $\left[ \Qm_\alpha, \h_\beta\right] =
\delta_{\alpha\beta}$.  This is the pairing that makes
$\Qsym$ and $\NCL$ graded dual Hopf algebras \cite{MR}.

The ribbon quasi-symmetric functions, $\Qs_\alpha$
are then defined
as the elements of $\Qsym$ such that
$\left[ \Qs_\alpha, \s_\beta\right] = \delta_{\alpha\beta}$.
Clearly, any $A \in \Qsym$ may be expanded in these bases
by using the formula
\begin{equation}
A = \sum_{\beta \models deg(A)} \left[ A, \s_\beta\right] \Qs_\beta
= \sum_{\beta \models deg(A)} \left[A, \h_\beta\right] \Qm_\beta.
\end{equation}
Similarly, for any element $\A \in \NCL$, $\A$ may be expanded
in terms of the $\s_\beta$ and $\h_\beta$ bases if we know
the pairing between $\A$ and $\Qs_\beta$ or $\Qm_\beta$.
\begin{equation}
\A = \sum_{\beta \models deg(\A)} \left[ \Qs_\beta, \A \right]
\s_\beta
= \sum_{\beta \models deg(\A)} \left[\Qm_\beta, \A\right] 
\h_\beta
\end{equation}

There is a simple relation between the $\Qsym/\NCL$ pairing
and the scalar product on $\NCL$.  This is expressed
with the following proposition.

\begin{prop}  \label{dualbase}
$A \in \Qsym$ and $\A \in \NCL$ are such that
for all $\B \in \NCL$ 
\begin{equation*}
\left[ A, \B \right] = \left< \A, \B \right>, \label{seq}
\end{equation*}
if and only if $A$ and $\A$ have the following relationship
\begin{equation*}\A = \sum_{\beta \models n}
(-1)^{\ell(\beta)+1} 
\left[ A, \s_\beta \right] \s_{\invc{\beta}} =
\sum_{\beta \models n}
(-1)^{\ell(\beta)+1}
\left[ A, \h_\beta \right] \h_{\invc{\beta}},
\end{equation*}
\begin{equation}
A  = \sum_{\beta \models n}
\left< \A, \s_\beta \right> \Qs_\beta =
\sum_{\beta \models n}
\left< \A, \h_\beta \right> \Qm_\beta.\label{makedual}
\end{equation}
\end{prop}

\begin{proof}  By applying Corollary \ref{scalarexpr},
\begin{align}
\left< \A, \B \right> &=
\left< \sum_{\beta \models n}
(-1)^{\ell(\beta)+1} 
\left[ A, \s_\beta \right] \s_{\invc{\beta}}, \B \right> \\
\nonumber
&= \left[ A, \sum_{\beta \models n}
(-1)^{\ell(\beta)+1}  
\left<\s_{\invc{\beta}}, \B \right> \s_\beta \right]
= \left[ A, \B \right].
\end{align}
\end{proof}

This last proposition implies that for a basis $\A_\alpha$
such that $\left<\A_\alpha,\A_\beta \right> = 
(-1)^{n+\ell(\alpha)} \delta_{\alpha\invc{\beta}}$, then
to compute its dual basis in $\Qsym$ it is sufficient to
calculate the values of the scalar products $\left< \A_\alpha, \s_\beta
\right>$, since by equation (\ref{makedual}), we have
that $\left[(-1)^{n + \ell(\invc{\alpha})} 
 \sum_{\beta \models n}
\left< \A_{\invc{\alpha}}, \s_\beta \right>\Qs_\beta, \A_\beta\right] =
(-1)^{n + \ell(\invc{\alpha})} 
\left< \A_{\invc{\alpha}}, \A_\beta \right> = 
\delta_{\alpha\beta}$.

\section{$q$-Analog bases}

\subsection{Scrambled Hopf algebra operators}
Consider the following transformation
on $Hom(A, A)$ that seems to arise in a natural
way when considered as an operation on symmetric functions
\cite{Z2}.  If $A$ is a Hopf algebra,
then for $V \in Hom(A,A)$ we define ${\ob{V}} = 
\mu \circ id \otimes (V\antipode) \circ \Delta = 
id \ast (VS)$.  We may show
that in any co-commutative Hopf algebra $A$, 
the bar operation on $V \in Hom(A,A)$ is an involution.
That is, we have

\begin{prop}  For any $V \in Hom(A,A)$ with $(A, \mu,
\unit,\Delta,\counit, \antipode)$ a co-com-mutative
Hopf algebra, we have that
$
{\ob{\ob{V}}} = V
$.
\end{prop}

Our $q$-analog arises by starting with a graded 
co-commutative Hopf algebra with a function $R^q$
such that for any element $f \in A$ that is of
homogeneous degree $R^q(f) = q^{deg(f)} f$.
The important properties of this function are that
$R^q\coeff_{q=0} = \unit\counit$ and $R^q \coeff_{q=1} = id$.
Now for any $V \in Hom(A,A)$, we set
\begin{equation}\label{qhatdef}
\hatq{V} = {\ob{\ob{V}~ \ob{R^q}}}.
\end{equation}

It is not at all obvious
that the effect on $V$ will be necessarily interesting.
We do however have many examples of where this $q$-analog
arises within the theory of symmetric functions \cite{Z}.  One
may only hope that within other co-commutative Hopf
algebras that this operation is also interesting.  In
this exposition we will show that there is an
application of this $q$-analog within the non-commutative
symmetric functions.

\begin{prop} \label{qanalp}
Let $A$ be a Hopf algebra, with the property 
that $\ob{V} := id \ast V\antipode$ 
is an involution for all $V \in Hom(A,A)$. 
Also assume that there is an operator $R^q$ 
such
that $R^q \coeff_{q=0} = \unit \counit$ and 
$R^q \coeff_{q=1} = id$.
If $V$ does not depend on $q$,
then the operator $\hatq{V} = \ob{\ob{V}~\ob{R^q}}$ 
is a $q$-analog of 
$V$ (recovered by setting $q=0$) and $\ob{\ob{V}~\unit 
\counit}$ (by setting $q=1$).
\end{prop}

\begin{proof}
Since $id \ast \antipode
= \unit \counit$, we have that 
$\ob{id} = \eta\epsilon$.  It
follows then that
\begin{equation}
\hatq{V} \coeff_{q=0} = 
\ob{\ob{V}~\ob{\unit\counit}} = \ob{\ob{V}} = V,
\end{equation}
\begin{equation}
\hatq{V} \coeff_{q=1} = 
\ob{\ob{V}~\ob{id}} = \ob{\ob{V}~\unit\counit}.
\end{equation}
\end{proof}

\subsection{A $q$-twisting of $\NCL$ operators}

Define an operator $A_\alpha$ that sends the ribbon
Schur function $\s_\beta$ to the ribbon Schur function
$\s_{\beta \append \alpha}$ (the concatenate operation)
and then extend this operator
linearly to the space $\NCL$.  Similarly define an
operator $B_\alpha$ that sends the ribbon Schur function
$\s_\beta$ to the ribbon Schur function $\s_{\beta \attach
\alpha}$ (the attach operation) 
and extend the function linearly.  

Define $\AAq$ as the $q$-analog of the operator 
$A_\alpha$ using equation
(\ref{qhatdef}).  From Proposition \ref{qanalp}, 
$\AAq$ has the property that
when $q=0$ it is $A_\alpha$.  When $q=1$ we see that
\begin{equation}
\AAq(\s_\beta)\coeff_{q=1} =
\ob{\ob{A_\alpha} \ncunit\nccounit}(\s_\beta) =
id \ast (\ob{A_\alpha} \ncunit\nccounit \ncantipode) 
(\s_\beta)
= \s_\beta \ob{A_\alpha}(1).
\end{equation}

It is easily shown that $\ob{A_\alpha}(1) = \s_\alpha$, 
hence we see that $\AAq(\s_\beta)\coeff_{q=1} = \s_\beta
\s_\alpha = A_\alpha(\s_\beta) + B_\alpha(\s_\beta)$.  In fact we may derive that there is an
simple formula for the action of $\AAq$.

\begin{prop} \label{Hqaction} Let $\beta$ be a 
composition of $n>0$, then
\begin{equation*}\AAq(\s_\beta) = A_\alpha(\s_\beta)
+ q^{n} B_\alpha(\s_\beta).
\end{equation*}
Also $\AAq(1) = \s_\alpha$.
\end{prop}

\begin{proof}
Define the numbers $\mC{\alpha}{\beta\gamma}$ as
the coefficients that arise in the coproduct $\ncDelta(\s_\alpha) =
\sum_{\beta,\gamma} \mC{\alpha}{\beta\gamma} \s_\beta\otimes
\s_\gamma$.  From the defining property of the antipode
we have the relation
\begin{equation}
\sum_{\beta\gamma} (-1)^{|\beta|} \mC{\alpha}{\beta\gamma}
\s_{\invb{\beta}} \otimes \s_{\gamma} = 
\sum_{\beta\gamma} (-1)^{|\gamma|} \mC{\alpha}{\beta\gamma}
\s_{\beta} \otimes \s_{\invb{\gamma}}=0\label{apodesum}
\end{equation}
as long as $\alpha$
is not the empty composition.

The notation for the $q$ analog of any operator is given
in equation (\ref{qhatdef}).  To show that the $q$-analog
of $A_\alpha$ satisfies the proposition we must give a
definition in terms of Hopf algebra operations and then
demonstrate that the action reduces dramatically.  For any
operator, (\ref{qhatdef}) may be restated as
\begin{equation}
\hatq{V} = \ncmu ( 1 \otimes V)
(\ncmu \otimes \ncantipode)( 1 \otimes \ncDelta)
( 1 \otimes \ncmu)
( 1 \otimes \ncantipode \otimes R^q)
( 1 \otimes \ncDelta ) \ncDelta.
\end{equation}


At this point it is a direct computation within the
Hopf algebra of the non-commutative symmetric functions
using relation (\ref{apodesum}) and the definition of
$A_\alpha$
to arrive at the formula stated in the proposition.
Since the computation is detailed and
not necessary for the remainder
of this exposition, we leave it to the reader as an
exercise.
\end{proof}

\subsection{An example of $q$-non-commutative symmetric
functions}
Define for a pair of compositions $\alpha,\beta \models n$
the statistic 
$c(\alpha,\beta) = \sum_{i \in D(\alpha) \intersect
D(\beta)} i$.

\vskip .1in

\begin{prop}\label{Haqdef}
Let $\H_\alpha^q = \sum_{\beta} q^{c(\alpha,\invc{\beta})} \s_\beta$
where the sum is over all compositions $\beta$ of 
$|\alpha|$ such that
$\alpha \finereq \beta$ then
$$\H_\alpha^q = \Axq{\alpha_{\ell(\alpha)}} 
\Axq{\alpha_{\ell(\alpha)-1}}
\cdots \Axq{\alpha_{1}} 1.$$
\end{prop}

Notice that when $q=0$, then $\H_\alpha^0 = \s_\alpha$.
We also have that when $q=1$, then $\H_\alpha^1 = 
\h_\alpha$.  We will think of this family as a
non-commutative analog of the Hall-Littlewood symmetric
functions because of the following two properties.  The
first says that this family shares a sort of self-duality
property similar to the $\s$, $\h$ and $\e$ bases of $\NCL$,
the second says that the commutative versions
of these symmetric functions agree with the Hall-Littlewood
symmetric functions when indexed by a composition that
is equivalent to a partition.

\begin{prop} \label{Hscalar}
\begin{equation}\left< \H_\alpha^q, \H_{\beta}^q \right>
= (-1)^{|\alpha|+\ell(\alpha)} \delta_{\alpha\invc{\beta}}.
\label{HaqHbq}\end{equation}
\end{prop}

\begin{proof}
Relation (\ref{HaqHbq}) will follow from Proposition \ref{Hqtscalar}
by setting $t=q$ and $q=0$.
\end{proof}

We also have the following remarkable connection with
the Hall-Littlewood basis, $H_{\lambda}^q$.

\begin{prop} 
If $\alpha = (1^a,b)$, then \label{Hqcomm}
$$\chi(\H_\alpha^q) = H_{(b,1^a)}^q.$$
\end{prop}

\begin{proof}
Due to the rule given in (\cite{M}, (5.7), p. 228)
we have the following important recurrence
for the Hall-Littlewood symmetric functions indexed
by a hook.
\begin{equation}
h_b H_{(1^a)}^q = H_{(b,1^a)}^q 
+ (1-q^{a}) H_{(b+1,1^{a-1})}^q
\end{equation}
Also consider the recurrence that we have developed
for these non-commutative symmetric functions.
\begin{align}
\H_{(1^a)}^q \h_b &= 
(A_{(b)} + q^a B_{(b)} +(1- q^a) B_{(b)})(\H_{(1^a)}^q) 
\nonumber\\
&= \H_{(1^a,b)}^q + (1-q^a) B_{(b)}(\H_{(1^a)}^q)
\end{align}
We also have that $\H_{(1^a)}^q = (A_{(1)} + q^{a-1} B_{(1)})
(\H_{(1^{a-1})}^q)$.  Since we have that $B_{(r)} B_{(s)} =
B_{(r+s)}$ and $B_{(r)} A_{(s)} = A_{(r+s)}$, then we see
$B_{(b)} (\H_{(1^a)}^q) = \H_{(1^{a-1},b+1)}^q$.

This implies
\begin{equation}
 \H_{(1^a)}^q \h_b = \H_{(1^a,b)}^q + 
(1-q^a) \H_{(1^{a-1},b+1)}^q.
\end{equation}
By induction on the length of the hook we see that
$\chi(\H_\alpha^q) = H_{(b,1^a)}^q$
(which is obviously
true when either $a=0$ or $b=1$).
\end{proof}

Using Proposition \ref{Hscalar} we are able to derive an equation
for the basis in $\Qsym$ dual to the $\H_\alpha^q$.  We notice
that $\delta_{\alpha\beta} = (-1)^{\ell(\alpha)+1}
\left< \H_{\invc{\alpha}}^q, \H_\beta^q \right>$.
This implies that if
\begin{equation}
\left[ \QPq_\alpha^q, \H_\beta^q \right] =
\delta_{\alpha\beta},
\end{equation}
then $\QPq_\alpha^q = \sum_{\beta\models |\alpha|}
\left<(-1)^{\ell(\alpha)+1}\H_{\invc{\alpha}}^q,
\s_\beta \right> \Qs_\beta$.  A simple calculation
using Proposition \ref{dualbase} implies
that
\begin{cor}
\begin{equation}
\QPq_\alpha^q = \sum_{\beta \finereq\alpha}
(-1)^{\ell(\beta)-\ell(\alpha)}
q^{c(\invc{\alpha},\beta)}
\Qs_\beta
\end{equation}
is the basis of $\Qsym$ which is dual to the family
$\H_\alpha^{q}$ with respect to the $\Qsym/\NCL$ pairing.
\end{cor}

We remark that these non-commutative symmetric 
functions are not equivalent to those defined in \cite{Hi}.
They are however remarkably similar and do agree for
$\alpha = (1^a,b)$.  We show this in the following
proposition.

Say that $D(\alpha) = 
\{ a_1 < a_2 < \ldots < a_{\ell(\alpha)-1}\}$
and $D(\beta) = 
\{ b_1 < b_2 < \ldots < b_{\ell(\beta)-1}\}$.
Let $Bre(\alpha,\beta)$ be the composition
of $\ell(\alpha)$ with the descent set equal to
$D(Bre(\alpha,\beta)) = \{\#\{a_j:
a_j \leq b_i\}: 1 \leq i \leq \ell(\beta)-1\}$.
Let $$\W_\alpha^q = \sum_{\beta \coarsereq \alpha}
q^{n(\invc{Bre(\alpha,\beta)} )} \s_\beta.$$
This is the definition of the non-commutative analogs
of the Hall-Littlewood symmetric functions given in
Theorem 6.13 of \cite{Hi}.  We may easily see that the
family $\W_\alpha^q$ satisfies the following recurrence

\begin{prop}
\begin{equation}
\W_{\alpha\append(m)}^q = A_{(m)}(\W_\alpha^q)
+q^{\ell(\alpha)} B_{(m)}(\W_\alpha^q)
\end{equation}
\end{prop}
\begin{proof}
$D(Bre(\alpha\append(m),\beta\append(m))) = D(Bre(\alpha,
\beta))\cup\{|\alpha|\}$.  At the same time we have
$D(Bre(\alpha\append(m),\beta\attach(m))) = D(Bre(\alpha,
\beta))$.  Both $Bre(\alpha\append(m),\beta\append(m))$
and $Bre(\alpha\append(m),\beta\attach(m))$ are compositions
of $\ell(\alpha)+1$, hence the proposition follows.
\end{proof}

This proposition should be compared to Proposition
\ref{Hqaction}. This also shows the following corollary.

\begin{cor}
For $\alpha = (1^a,b)$ we have
$$\W_\alpha^q = \H_\alpha^q.$$
\end{cor}

One open question that arises from this definition is: Is
there a Hecke algebra action on $\Qsym$ such that 
the functions $\QPq_\alpha^q$ are invariants
under its action?  This is the case of the functions
of Hivert that are dual to the non-commutative symmetric
functions $\W_\alpha^q$ and are given by the formula:

\begin{equation}
\QGq_{\alpha}^q =\sum_{\beta \finereq \alpha}
(-1)^{\ell(\beta) - \ell(\alpha)} q^{n(\invb{
Bre(\beta,\alpha)})}
\Qs_\beta.
\end{equation}

The functions $\H_\alpha^q$ have a factorization property
that is very similar to that held by the functions of
Hivert $\W_\alpha^q$ and the commutative Hall-Littlewood
symmetric functions \cite{Hi}.

\begin{thm}\label{factor}
Let $\zeta$ be an $r^{th}$ root of unity.  Then
$$\H_\alpha^{\zeta} = \H_{\alpha^{(1)}}^\zeta
\H_{\alpha^{(2)}}^\zeta \cdots \H_{\alpha^{(k)}}^\zeta$$
for any decomposition $\alpha = \alpha^{(1)}
\append \alpha^{(2)} \append \cdots \append \alpha^{(k)}$,
where for $1 \leq i \leq k-1$, $\alpha^{(i)}$ is a composition
of a multiple of $r$.
\end{thm}

This theorem follows from the following derivation
of the product rule for these functions and then evaluating
$q$ at a root of unity.

\begin{prop}\label{prodHq}
\begin{equation}
\H_\alpha^{q} \H_\beta^q = \sum_{\gamma \coarsereq \beta}
f_{\alpha\beta}^{\gamma}(q) (\H_{\alpha \append\gamma}^q
+(1-q^{|\alpha|}) \H_{\alpha\attach\gamma}^q),
\end{equation}
where
\begin{equation}
f_{\alpha\beta}^{\gamma}(q) = q^{c(\beta,\invc{\gamma})}
(1-q^{|\alpha|})^{\ell(\beta) - \ell(\gamma)}.
\end{equation}
\end{prop}

Before proceeding with the proof we introduce the following
lemma.

\begin{lem}\label{Hasb}
\begin{equation}
\H_\alpha^q \s_{\beta} = \sum_{\gamma \coarsereq \beta}
g_{\alpha\beta}^{\gamma}(q)(\H_{\alpha\append\gamma}^q
+(1-q^{|\alpha|})\H_{\alpha\attach\gamma}^q),
\end{equation}
where
\begin{equation}
g_{\alpha\beta}^{\gamma}(q) = (-1)^{\ell(\beta)-\ell(\gamma)}
q^{c(\alpha\append\beta,\invc{(\alpha\append\gamma)})}.
\end{equation}
\end{lem}

\begin{proof}
Let $n=|\alpha|+|\beta|$.
We will take the scalar product defined in equation
(\ref{NCscalar}) of $\H_{\alpha}^q \s_\beta$ and
$\H_{\invc{\theta}}^q$.  This will give the coefficient of
$ (-1)^{\ell(\theta)+n} \H_{\theta}^q$ in the
expression $\H_{\alpha}^q \s_\beta$.  By expanding $\H_\alpha^q$
in terms of $\s_\gamma$ and using that 
$\left< \s_\alpha, \H_\beta^q \right>
= (-1)^{|\alpha|+\ell(\alpha)} q^{c(\alpha,\beta)} 
\chi(\alpha \finereq
\invc{\beta})$, we see
\begin{align}
(-1)^{\ell(\theta)+n}\left< \H_{\alpha}^q \s_\beta, 
\H_{\invc{\theta}}^q \right> &=
(-1)^{\ell(\theta)+n}
\sum_{\gamma\coarsereq\alpha} q^{c(\alpha,\invc{\gamma})}
\left< \s_{\gamma\append\beta} + \s_{\gamma\attach\beta},
\H_{\invc{\theta}}^q \right> \nonumber \\
&=
\sum_{\gamma\coarsereq\alpha} q^{c(\alpha,\invc{\gamma})}
(-1)^{\ell(\theta)+\ell(\gamma\append\beta)} 
q^{c(\gamma\append\beta,\invc{\theta})} 
\chi(\gamma\append\beta \finereq \theta)\\
&\hskip .1in+ \sum_{\gamma\coarsereq\alpha} 
q^{c(\alpha,\invc{\gamma})}
(-1)^{\ell(\theta)+\ell(\gamma\attach\beta)} 
q^{c(\gamma\attach\beta,\invc{\theta})} 
\chi(\gamma\attach\beta \finereq \theta).\nonumber
\end{align}

Break the composition $\theta$ into the composition consisting
of the first $|\alpha|$ `cells' of $\theta$: $\mu = \theta
\coeff_{1\ldots|\alpha|}$, and the last $|\beta|$ `cells:'
$\nu = \theta\coeff_{|\alpha|+1\ldots n}$, so that
either $\theta = \mu \append\nu$
or $\theta =  \mu \attach\nu$.  

If $\nu \finer \beta$ or
$\mu \finer \alpha$ then clearly both sums will be $0$.

If $\mu \coarser \alpha$ then $D(\invc{\mu}) \cap D(\alpha)$
is non empty and contains at least one element $a$.  
For each $\gamma$ in the sums, either $D(\gamma)$ contains
$a$ or it doesn't.  The terms with $\gamma$ such that
$a \in D(\gamma)$ are of opposite sign but the same $q$
coefficient as those $\gamma$ such that $a \notin D(\gamma)$.
Therefore the two sums will again be $0$.

We need only consider the $\theta$ where $\mu = \alpha$ 
and $\nu \coarsereq \beta$.  If $\theta = \alpha\append\nu$,
then the second sum is clearly $0$ and first sum contains
only $1$ term, $(-1)^{\ell(\nu)+\ell(\beta)}
q^{c(\alpha\append\beta,\invc{\theta})}$ (which agrees
with the statement of the lemma).
If $\theta = \alpha\attach\nu$, then both sums have exactly
one non-zero term, and the scalar product is
\begin{align}
&(-1)^{\ell(\nu)+\ell(\beta)-1} 
q^{c(\alpha\append\beta,\invc{(\alpha\attach\nu)})}+
(-1)^{\ell(\nu)+\ell(\beta)} 
q^{c(\alpha\attach\beta,\invc{(\alpha\attach\nu)})}\\
&\hskip .2in= (-1)^{\ell(\nu)+\ell(\beta)}
q^{c(\alpha\attach\beta,\invc{(\alpha\attach\nu)})}
(1-q^{|\alpha|}).\nonumber
\end{align}
\end{proof}

\begin{proof}[Proof of Proposition \ref{prodHq}] 
Expanding $\H_\beta^q$ in terms of $\s_\gamma$ and using
Lemma \ref{Hasb}, yields
\begin{align}\nonumber
&(-1)^{\ell(\theta)+n}\left< \H_{\alpha}^q \H_\beta^q, 
\H_{\invc{\theta}}^q \right>\\
&=(-1)^{\ell(\theta)+n}
\sum_{\gamma \coarsereq \beta} \sum_{\mu \coarsereq \gamma}
q^{c(\beta,\invc{\gamma})}
\left<g_{\alpha\gamma}^{\mu}(q) 
(\H_{\alpha \append\mu}^q
+(1-q^{|\alpha|}) \H_{\alpha\attach\mu}^q),
\H_{\invc{\theta}}^q \right>\\
&=
\sum_{\gamma \coarsereq \beta} \sum_{\mu \coarsereq \gamma}
q^{c(\beta,\invc{\gamma})}
g_{\alpha\gamma}^{\mu}(q) 
\delta_{\alpha \append\mu,\theta}
+q^{c(\beta,\invc{\gamma})}
g_{\alpha\gamma}^{\mu}(q) (1-q^{|\alpha|})
\delta_{\alpha\attach\mu,\theta}.\nonumber
\end{align}
Now if $|\alpha|$ is in $D(\theta)$ then  the inner
product is
\begin{align}\nonumber
&=\sum_{\mu \coarsereq \gamma \coarsereq \beta}
q^{c(\beta,\invc{\gamma})}
(-1)^{\ell(\gamma) -\ell(\mu)}
q^{c((\alpha\append\gamma),\invc{(\alpha\append\mu)})}\\
&=
\sum_{\mu \coarsereq \gamma \coarsereq \beta}
q^{c(\beta,\invc{\gamma})+c(\gamma,\invc{\mu})}
(-q^{|\alpha|})^{\ell(\gamma) -\ell(\mu)}\\
&=q^{c(\beta,\invc{\mu})} (1-q^{|\alpha|})^{\ell(\beta)-
\ell(\mu)}.\nonumber
\end{align}
This agrees with the formula given for 
$f_{\alpha\beta}^\mu(q)$.  If $|\alpha|$ is 
not in $D(\theta)$, the result is $(1-q^{|\alpha|})$
times this result.
\end{proof}

\section{$q,t$-Analogs of non-commutative symmetric functions}
Define the following $q,t$-non commutative symmetric function.
\begin{equation}
\H_{\alpha}^{qt} = \sum_{\beta\models|\alpha|}
t^{c(\alpha,\invc{\beta})} q^{c(\inva{\alpha},\invb{\beta})}
\s_\beta
\end{equation}

Clearly from this definition, if $q=0$ and $t=q$, 
then all terms such that $D(\inva{\alpha}) \cap D(\invb{\beta})$
is non empty vanish and 
we have $\H_{\alpha}^{0q} = \H_\alpha^q$.
Therefore we
also have the specializations, $\H_{\alpha}^{00} = \s_\alpha$
and $\H_{\alpha}^{01} = \h_\alpha$, and $\H_\alpha^{10} 
= \e_{\invc{\alpha}}$.

Moreover, $\H_\alpha^{qt}$ satisfies the following relations
which are similar to those held by the Macdonald symmetric
functions in the commutative case:
\begin{equation}
\H_\alpha^{tq} = \oma \H_{\inva{\alpha}}^{qt}.
\end{equation}
\begin{equation}
q^{n(\inva{\alpha})} t^{n(\alpha)} 
\H_\alpha^{\frac{1}{q}\frac{1}{t}} =
\omc \H_{\alpha}^{qt}.
\end{equation}

When we set $q=1$, the $\H_{\alpha}^{1t}$ become products
of some non-standard non-commutative symmetric functions, 
as seen in the following proposition.

\begin{prop} \label{H1tprod}
Define the non-commutative symmetric functions
$\H_{(m)}^{q(i)} = \sum_{\beta \models m} 
q^{(\ell(\beta)-1) i
+ n(\invb{\beta})} \s_\beta$. For a composition $\alpha$ such that
$k = \ell(\alpha)$, we have
\begin{equation}
\H_{\alpha}^{q1} =
\H_{(\alpha_1)}^{q(\sum_{i>1}\alpha_i)}
\H_{(\alpha_{2})}^{q(\sum_{i>2}\alpha_i)}
\cdots
\H_{(\alpha_k)}^{q(0)}.\label{factorize}
\end{equation}
\end{prop}

\begin{proof}
Fix $\alpha$ and for $1 \leq i \leq \ell(\alpha)$
let $\gamma^{(i)}$ be a composition of $\alpha_i$.
The coefficient of $\s_{\gamma^{(1)}} 
\s_{\gamma^{(2)}} \cdots \s_{\gamma^{(\ell(\alpha))}}$
in the right hand side of equation (\ref{factorize})
is $q$ raised to the power of
\begin{equation}
\sum_i n(\invb{\gamma^{(i)}}) + 
\sum_i (\ell(\gamma^{(i)})-1)\sum_{j>i} \alpha_i 
= \sum_i \left(\sum_{d \in D(\invb{\gamma^{(i)}})} i + 
\sum_{j>i} \alpha_i\right).
\end{equation}
This agrees with $c(\inva{\alpha},\invb{\beta})$ where
$\beta$ is attach and concatenate of the $\gamma^{(i)}$
and hence agrees with the $q$ coefficient on the left
hand side of equation (\ref{factorize}).
\end{proof}

We also have the following two additional Propositions
that lead us to believe that it is an interesting
generalization of the family $\H_\alpha^q$.

\begin{prop} Let $\alpha \models n$, then\label{Hqtscalar}
\begin{equation}
\left< \H_{\alpha}^{qt}, \H_{\beta}^{qt} \right>
= (-1)^{|\alpha| + \ell(\alpha)}
\delta_{\alpha\invc{\beta}}
\prod_{i=1}^{n-1} (1-q^i t^{n-i}) .
\end{equation}
\end{prop}

\begin{proof}
\begin{align}
\left< \H_{\alpha}^{qt}, \H_{\beta}^{qt} \right>
&= \sum_{\gamma} \sum_{\theta} t^{c(\alpha,\invc{\gamma})+
c(\beta,\invc{\theta})} q^{c(\inva{\alpha},\invb{\gamma})+
c(\inva{\beta}, \invb{\theta})} \left< \s_\gamma, 
\s_\theta \right> \\
&= \sum_{\gamma} (-1)^{n+\ell(\gamma)}
t^{c(\alpha,\invc{\gamma})+
c(\beta,\gamma)} q^{c(\inva{\alpha},\invb{\gamma})+
c(\inva{\beta}, \inva{\gamma})}. \nonumber
\end{align}
If $\alpha \neq \invc{\beta}$, then $(D(\alpha) \cap D(\beta))
\cup (D(\invc{\alpha}) \cap D(\invc{\beta}))$
is non empty.  Take the smallest element $i$ of this set 
(although any will do) and consider the involution $\phi$ on the
set of compositions such that the compositions
that contain $i$ in the
descent set are sent to the compositions that do not
contain the element $i$ (in the most natural manner). 
For each $\gamma \models n$, the terms corresponding to
$\gamma$ and $\phi(\gamma)$ have the same weight but
opposite sign, hence the sum is $0$.

If $\alpha = \invc{\beta}$, then the sum reduces to
\begin{align}
&=\sum_{\gamma} (-1)^{n+\ell(\gamma)}
t^{c(\alpha,\invc{\gamma})+
c(\invc{\alpha},\gamma)} q^{c(\inva{\alpha},\invb{\gamma})+
c(\invb{\alpha}, \inva{\gamma})} \\
&= \sum_{S \subseteq \{ 1,\ldots,n\}} 
(-1)^{n+|S \cap D(\invc{\alpha})| + |S^c \cap D(\alpha)| + 1} 
t^{\sum_{i \in S} i} q^{\sum_{i \in S} n+1-i},\nonumber
\end{align}
where the subsets $S$ represent the sets $(D(\alpha) \cap
D(\invc{\gamma})) \cup (D(\invc{\alpha}) \cap D(\gamma))$.
This is clearly equal to the product stated in the 
proposition.
\end{proof}

\begin{cor}  The family
\begin{equation}
\QPq_\alpha^{qt} = \prod_{i=1}^{n-1}\frac{1}{1-q^i t^{n-i}}
\sum_{\beta}
(-1)^{\ell(\beta)-\ell(\alpha)}
t^{c(\invc{\alpha},\beta)} q^{c(\invb{\alpha},\inva{\beta})}
\Qs_\beta,
\end{equation}
has the property that $\left[ \QPq_{\alpha}^{qt},
\H_\beta^{qt} \right] = \delta_{\alpha\beta}$.
\end{cor}
\begin{proof}
If we wish that $\left[ \QPq_{\alpha}^{qt},
\H_\beta^{qt} \right] = \delta_{\alpha\beta}$, then
using equation (\ref{makedual}), 
$$\QPq_\alpha^{qt} 
= \sum_{\beta \models 
|\alpha|} \left< \A, \s_\beta \right> \Qs_\beta,$$
where
$\A = (-1)^{\ell(\alpha)+1} \prod_{i=1}^{n-1}
\frac{1}{1-q^it^{n-i}}
\H_{\invc{\alpha}}^{qt}$, since
\begin{equation}
\delta_{\alpha\beta} = (-1)^{\ell(\alpha)+1} 
\prod_{i=1}^{n-1}
\frac{1}{1-q^it^{n-i}}\left<
\H_{\invc{\alpha}}^{qt},
\H_\beta^{qt}\right>.
\end{equation}
A simple calculation yields the equation stated in the
corollary.
\end{proof}

\def\P{\hbox{\scr P}}

There is a characterization of the non-commutative
$q,t$ analogs $\H_\alpha^{qt}$ in terms of 
properties that are similar to those shared by the
Macdonald symmetric functions.  This characterization
is not particularly important for our treatment, but
it should not be ignored because of the similarities
that it shares with the commutative case.

Define a family of non-commutative symmetric functions
$\P_\alpha^{qt}$ by the following three conditions.

\vskip .2in
\noindent
1. $\P_\alpha^{qt} = \H_\alpha^{t} + \sum_{\beta<\alpha}
c_{\alpha\beta}(q,t) \H_\beta^t$ for some coefficients 
$c_{\alpha\beta}(q,t)$ that are rational functions in
the parameters $q$ and $t$.

\vskip .2in
\noindent
2. $\oma \P_\alpha^{qt} = a_\alpha(q,t)
\P_{\inva{\alpha}}^{tq}$ for some coefficients
$a_\alpha(q,t)$.

\vskip .2in
\noindent
3. $\left< \P_\alpha^{qt}, \P_\beta^{qt} \right> = 0$
if $\alpha \neq \invc{\beta}$.
\vskip .2in
\begin{thm}
The family $\P_\alpha^{qt}$ are defined by the
three conditions listed above and, moreover, 
$\P_\alpha^{qt} = r_\alpha
\H_\alpha^{qt}$ where $r_\alpha = 1/\prod_{i \in
D(\invc{\alpha})}(1-q^{n-i}t^{i})$.  The coefficients
$c_{\alpha\beta}(q,t)$ are given by the formula
$$c_{\alpha\beta}(q,t)=\prod_{i
\in D(\invc{\alpha})\cap D(\beta)}q^{n-i}/(1-t^i
q^{n-i}).$$  The coefficients
$a_\alpha(q,t)$ mentioned in the second condition are given
by the formula
$a_\alpha(q,t) = \prod_{i \in D(\alpha)} (1-q^{n-i}
t^i)/\prod_{i \in D(\invc{\alpha})} (1-q^{n-i}t^i)$.
\end{thm}

\begin{proof}
The proof proceeds by induction, for there is a method of
calculating the coefficients $c_{\alpha\beta}(q,t)$ from
ones preceding it in some order.  Say that $\P_\alpha^{qt}
= \sum_{\beta \leq \alpha} c_{\alpha\beta}(q,t) \H_\beta^t$
where we assume that $c_{\alpha\alpha}(q,t) = 1$ and that
this family satisfies the three conditions given above.

Assume that the coefficients $c_{\gamma\delta}(q,t)$ are
known and given by the formula stated in the theorem for
all $\gamma$ such that $|D(\gamma)|>|D(\alpha)|$ or for
$\gamma=\alpha$ and $\delta > \beta$.  To determine
$c_{\alpha\beta}(q,t)$ we take the scalar product
of $\P_\alpha^{qt}$ and $\oma \P_{\invb{\beta}}^{tq}$
since $\beta < \alpha$, $|D(\invb{\beta})|>|D(\alpha)|$
and all coefficients in
$\P_{\invb{\beta}}^{qt}$ have been calculated already.
Since $\oma \P_{\invb{\beta}}^{tq} = \P_{\invc{\beta}}^{qt}
= \sum_{\theta \leq \invc{\beta}} c_{\invc{\beta}\theta}(q,t)
\H_\theta^t$, hence we have the expression 

\begin{align*}
\left<\P_{\alpha}^{qt},\oma \P_{\invb{\beta}}^{tq} 
\right> &= \left< \P_\alpha^{qt},
a_{\invb{\beta}}(t,q) \P_{\invc{\beta}}^{qt}
\right> = 0\\ &= 
a_{\invb{\beta}}(t,q)
\sum_{\alpha \geq \theta
\geq \beta} c_{\alpha\theta}(q,t)
c_{\invc{\beta}\invc{\theta}}(q,t)
(-1)^{n+\ell(\theta)}\\
&=a_{\invb{\beta}}(t,q)c_{\alpha\beta}(q,t)
(-1)^{n+\ell(\beta)}+\\
& \hskip .3ina_{\invb{\beta}}(t,q)
\sum_{\alpha \geq \theta
> \beta} c_{\alpha\theta}(q,t)
c_{\invc{\beta}\invc{\theta}}(q,t)
(-1)^{n+\ell(\theta)}.
\end{align*}

Those values of $c_{\invc{\beta}\theta}(q,t)$ may be
calculated from what we have already determined since
\begin{equation*}
c_{\invc{\beta}\invc{\theta}}(q,t) =
\left< \H_{\theta}^t, \P_{\invc{\beta}}^{qt}
\right>(-1)^{n+\ell(\theta)} =
\left< \H_{\theta}^t, \oma \P_{\invb{\beta}}^{tq}
\right>/a_{\invb{\beta}}(t,q)(-1)^{n+\ell(\theta)}.
\end{equation*}

Hence we see that
\begin{equation*}
c_{\alpha\beta}(q,t) =
\frac{(-1)^{n+1+\ell(\beta)}}{a_{\invb{\beta}}(t,q)}
\sum_{\alpha \geq \theta
> \beta} c_{\alpha\theta}(q,t)
\left< \H_{\theta}^t, \oma \P_{\invb{\beta}}^{tq}
\right>.
\end{equation*}
In addition we may calculate $a_{\invb{\beta}}(t,q)$
since
\begin{equation*}
\left< \H_\beta^{t}, \omega \P_{\invb{\beta}}^{tq} \right>
=
\left< \H_\beta^t, a_{\invb{\beta}}(t,q)
\P_{\invc{\beta}}^{qt} \right>
= a_{\invb{\beta}}(t,q) (-1)^{n + \ell(\beta)}.
\end{equation*}
 Although we may use these formulas to
calculate the coefficients, the only conclusion that we are
going to draw from them is that the coefficients
$c_{\alpha\beta}(q,t)$ are determined by assuming that the
defining conditions are true, hence the family
$\P_\alpha^{qt}$ which satisfies these conditions is unique.

It remains to show then that $\P_\alpha^{qt} =
\H_\alpha^{qt}/\prod_{i \in
D(\invc{\alpha})}(1-q^{n-i}t^{i})$ satisfies the
conditions listed above.  Clearly they satisfy
conditions 2 and 3.  It remains to show that the
$\H_\alpha^{qt}$ have the correct expansion in terms
of $\H_\alpha^t$.

\begin{align*}
\left<\H_{\invc{\alpha}}^t, \H_\alpha^{qt}\right>
&=\sum_{\beta \geq \invc{\alpha}} t^{c(\invc{\alpha},
\invc{\beta})}\left<\s_\beta, \H_\alpha^{qt}\right>\\
&=\sum_{\beta \geq \invc{\alpha}} (-1)^{n+\ell(\beta)}
t^{c(\alpha,\beta)+c(\invc{\alpha},\invc{\beta})}
q^{c(\inva{\alpha},\inva{\beta})}\\
&=\sum_{\beta \geq \invc{\alpha}} (-1)^{n+\ell(\beta)}
t^{c(\invc{\alpha},\invc{\beta})}
q^{c(\inva{\alpha},\inva{\beta})}\\
&= \sum_{S \subseteq D(\invc{\alpha})}
(-1)^{n+1+|S|}
t^{\sum_{i \in S} i}q^{ \sum_{i\in S} n-i}\\
&= (-1)^{n+\ell(\invc{\alpha})} \prod_{i \in
D(\invc{\alpha})}(1-q^{n-i}t^{i})
\end{align*}

We also see for $\beta$ is not strictly smaller than
$\alpha$ then $D(\alpha ) \cap D(\invc{\beta} )$
is non-empty and
\begin{align*}
\left<\H_{\invc{\beta}}^t, \H_\alpha^{qt}\right>
&=\sum_{\gamma \geq \invc{\beta}} t^{c(\invc{\beta},
\invc{\gamma})}\left<\s_\gamma, \H_\alpha^{qt}\right>\\
&=\sum_{\gamma \geq \invc{\beta}} t^{c(\invc{\beta},
\invc{\gamma})}
t^{c(\alpha,\gamma)}
q^{c(\inva{\alpha},\inva{\gamma})} (-1)^{n+\ell(\gamma)}.
\end{align*}
Since there is some element $a$ in $D(\alpha ) \cap
D(\invc{\beta} )$, every composition $\gamma \geq
\invc{\beta}$ 
either has $a \in D(\gamma)$ or $a \in D(\invc{\gamma})$. 
There is an obvious involution between these two sets of
compositions and they have the same weight but opposite
sign, hence the sum is $0$ in this case.
\end{proof}

As in the case of the family with one parameter, when
the functions are indexed by composition representing
a partition (i.e. a hook), 
then they are a generalization of
the Macdonald symmetric function.

\begin{prop}  If $\alpha = (1^a,b)$, then \label{Hqtcomm} 
\begin{equation}
\chi(\H_\alpha^{qt}) = H_{(b,1^a)}^{qt}.
\end{equation}
\end{prop}

\begin{proof}
Idea: same as in $q$ case.  Show that 
\begin{equation}
\H_{(1^a)}^{qt} \H_{(b)}^{qt} =
\frac{1-t^a}{1-q^b t^a} \H_{(1^{a-1},b+1)}^{qt}
+ \frac{1-q^b}{1-q^b t^a} \H_{(1^{a},b)}^{qt}
\end{equation}
and by a formula (\cite{M} eq. (6.24) p. 340) we have the same recurrence in
the commutative case.  That is,
\begin{equation}
H_{(1^a)}^{qt} H_{(b)}^{qt} = 
\frac{1-t^a}{1-q^b t^a} H_{(b+1,1^{a-1})}^{qt}
+ \frac{1-q^b}{1-q^b t^a} H_{(b,1^{a})}^{qt}.
\end{equation}
By induction this implies
that the commutative versions agree on hooks.

Consider the product $\H_{(1^a)}^{qt} \H_{(b)}^{qt}$. 
This is equal to
\begin{align}
\H_{(1^a)}^{qt} \H_{(b)}^{qt} &=
\sum_{\gamma \models a} \sum_{\theta \models b}
t^{n(\invc{\gamma})} 
q^{n(\invb{\theta})}
\s_{\gamma} \s_\theta \\
&= \sum_{\gamma \models a} \sum_{\theta \models b}
t^{n(\invc{\gamma})} 
q^{n(\invb{\theta})}
(A_\theta(\s_\gamma) + 
 B_{\theta}(\s_\gamma) ).\nonumber
\end{align}
We also have since $\inva{(1^a,b)} = (1^{b-1},a+1)$.
\begin{align}
\H_{(1^{a},b)}^{qt} &= \sum_{\beta\models a+b}
t^{c((1^a,b),\invc{\beta})} 
q^{c((1^{b-1},a+1),\invb{\beta})}
\s_\beta\\
&=\sum_{\gamma \models a} \sum_{\theta \models b}
t^{n(\invc{\gamma})} 
q^{n(\invb{\theta})}  (A_\theta(\s_\gamma) + 
t^{a} B_{\theta}(\s_\gamma) ).\nonumber
\end{align}
While at the same time
\begin{align}
\H_{(1^{a-1},b+1)}^{qt} &= \sum_{\beta\models a+b}
t^{c((1^{a-1},b+1),\invc{\beta})} 
q^{c((1^{b},a),\invb{\beta})}
\s_\beta\\
&=\sum_{\gamma \models a} \sum_{\theta \models b}
t^{n(\invc{\gamma})} 
q^{n(\invb{\theta})}  (q^{b} A_\theta(\s_\gamma) + 
B_{\theta}(\s_\gamma) ).\nonumber
\end{align}
From here is easily shown that
\begin{equation}
(1-q^b t^a )\H_{(1^a)}^{qt} \H_{(b)}^{qt} =
(1-q^b)\H_{(1^{a},b)}^{qt} +
(1-t^a) \H_{(1^{a-1},b+1)}^{qt}.
\end{equation}
\end{proof}

These two properties are only an indication that
$\H_{\alpha}^{qt}$ are an important generalization of
the Macdonald symmetric functions.  The
first property does not occur in many families of
non-commutative symmetric functions, the second, however,
could
appear for many different families (since the functions
of Hivert also have the same property that they have
the 'correct' statistic on hooks).

The most important indication that  the family 
$\H_{\alpha}^{qt}$ is
indeed an important analog to the Macdonald symmetric
functions is the appearance of an operator `nabla' that
is analogous to the operator introduced in
\cite{BG} for the symmetric functions.

First we define the analog $\Hw_\alpha^{qt} =
t^{n(\alpha)} \H_\alpha^{q\frac{1}{t}} =
\sum_{\beta\models|\alpha|}
t^{c(\alpha,\beta)} q^{c(\inva{\alpha},\invb{\beta})}
\s_\beta$.
Next, define $\nabla$ to be a linear operator
with the property that $\nabla (\Hw_\alpha^{qt}) =
t^{n(\alpha)}
q^{n(\inva{\alpha})} \Hw_\alpha^{qt}
$.  For our scalar product, we have that
\begin{equation} \label{nabscal}
\left< \nabla(f), \nabla(g) \right> = 
q^{n \choose {2}} t^{{n} \choose{2}} \left< f, g\right>.
\end{equation}

As we will see, this operator has many properties
that are analogous to those seen in the commutative case.  
In the non-commutative case the situation is somewhat
simpler and we are able to state precisely the action
of the operator on the ribbon basis.

\begin{prop}\label{nablaact} If $\alpha \models n$, then
\begin{equation}
\nabla(\s_\alpha) = 
(-1)^{n+\ell(\alpha)}
q^{n(\inva{\alpha})}
t^{n(\invc{\alpha})}
\sum_{\beta\finereq\invc{\alpha}} 
\prod_{i \in D(\alpha)\cap
D(\beta)} (t^i + q^{n-i}) \s_\beta. \label{nabfor}
\end{equation}
\end{prop}

To prove this formula we will need several lemmas for
the action of these operators on various bases.  By
choosing good notation for these operators, the proofs
become almost transparent.  We will order the ribbons
by their descent sets using the total order described
section \ref{comps}.

For a tensor product of two matrices with $B = \left[
b_{ij} \right]_{1 \leq i,j \leq n}$ we will use the
convention that
\begin{equation}
A \otimes B = \left[ b_{ij} A \right]_{1 \leq i,j \leq n}.
\end{equation}
That is, the $(r,s)$ entry in this matrix is
$$b_{(r~div~n)+1, (s~div~n)+1} a_{(r~mod~n)+1, (s~mod~n)+1}.$$

\begin{lem}
Let $\s$ be a column vector of $\s_\alpha$ and $\Hw$ is
a column vector of $\Hw_\alpha^{qt}$ (both using the
total order of section \ref{comps}), then
\begin{equation}
\Hw = \left[\begin{matrix} 1 & q^{n-1} \\
1 & t
\end{matrix}\right] \otimes  \left[\begin{matrix} 1 & q^{n-2} \\
1 & t^2
\end{matrix}\right] \otimes 
\cdots \otimes  \left[\begin{matrix} 1 & q \\
1 & t^{n-1}
\end{matrix}\right]\s.\label{Hwtos}
\end{equation}
\end{lem}

This lemma follows by realizing that if $\phi(\alpha) = k$,
then the entries in
the $k^{th}$ row of the tensor product matrix
agrees with the formula for the coefficients
of $\s_\beta$ in $\Hw_\alpha^{qt}$.  By taking the inverse
of this tensor product matrix we derive the inverse relation.

\begin{cor}
Let $\s$ be a column vector of $\s_\alpha$ and $\Hw$ is
a column vector of $\Hw_\alpha^{qt}$, then
\begin{equation}
\s = \left(\prod_{i=1}^{n-1} \frac{1}{t^i-q^{n-i}}\right) 
\left[\begin{matrix} t & -q^{n-1} \\
-1 & 1
\end{matrix}\right] \otimes 
\left[\begin{matrix} t^2 & -q^{n-2} \\
-1 & 1
\end{matrix}\right] \otimes 
\cdots \otimes  \left[\begin{matrix} t^{n-1} & -q \\
-1 & 1
\end{matrix}\right]\Hw.\label{stoHw}
\end{equation}
\end{cor}

\begin{lem}
Let $\Hw$ be a column vector of $\Hw_\alpha^{qt}$, then
\begin{equation}
\nabla \Hw = \left[\begin{matrix} q^{n-1} & 0 \\
0 & t
\end{matrix}\right] \otimes  \left[\begin{matrix} q^{n-2} & 0 \\
0 & t^2
\end{matrix}\right] \otimes 
\cdots \otimes  \left[\begin{matrix} q & 0 \\
0 & t^{n-1}
\end{matrix}\right]\Hw. \label{nablaonHw}
\end{equation}
\end{lem}

The proof of this lemma again follows by calculating the
entry in the row indexed by $\phi(\alpha)$.

\begin{proof}{(of Proposition \ref{nablaact})}
We calculate the action of $\nabla$ on the column
vector $\s$.  This follows by first expressing $\s$
in terms of $\Hw$ using equation (\ref{stoHw}), then
using the action of $\nabla$ on $\Hw$, then reexpressing
the answer in terms of $\s$ using (\ref{Hwtos}).  We
calculate that
\begin{equation}
\frac{1}{t^i-q^{n-i}}
 \left[\begin{matrix} t^{i} & -q^{n-i} \\
-1 & 1
\end{matrix}\right]
\left[\begin{matrix} q^{n-i} & 0 \\
0 & t^{i}
\end{matrix}\right]
\left[\begin{matrix} 1 & q^{n-i} \\
1 & t^{i}
\end{matrix}\right] =
\left[\begin{matrix} 0 & -q^{n-i} t^{i} \\
1 & (t^{i} + q^{n-i})
\end{matrix}\right].
\end{equation}
Therefore we see that the action of $\nabla$ on $\s$
is given by the equation
\begin{equation}
\nabla(\s) =
\left[\begin{matrix} 0 & -q^{n-1} t \\
1 & (t + q^{n-1})
\end{matrix}\right] \otimes
\left[\begin{matrix} 0 & -q^{n-2} t^2 \\
1 & (t^2 + q^{n-2})
\end{matrix}\right] \otimes
\cdots \otimes
\left[\begin{matrix} 0 & -q t^{n-1} \\
1 & (t^{n-1} + q)
\end{matrix}\right] \s.
\end{equation}

Now translate this tensor product directly to the action
on the $\s_\alpha$ basis to arrive at the formula stated in
equation (\ref{nabfor}).
\end{proof}

Interesting connections arise with combinatorics and 
representation theory that are analogous to the commutative
case.  Recall that for the standard symmetric functions
we have that $\left< \nabla(e_n), h_{1^n} \right>$
is a $q,t$ analog of the number $(n+1)^{n-1}$, which is the
number of parking functions (a function $f:[n] \rightarrow [n]$
is a parking function if the sequence 
$(f(1), f(2), \ldots, f(n))$
when sorted in increasing order $(a_1, a_2, \ldots, a_n)$
satisfies $a_i \leq i$). We also know 
that $\nabla(e_n)
\coeff_{t=1}$ is a graded Frobenius series for the parking
function module \cite{H}.  Moreover, the top component 
$\nabla(e_n) \coeff_{e_n}$ is known to be
a $q,t$-analog of the number of increasing parking functions
which is given by the Catalan numbers \cite{GH}.

In the non-commutative case, these statements occur with
exact analogy.  We will see that analogs of the parking 
functions are the preferential arrangements (\cite{St}
p. 146).    An exponential generating function for
the number of preferential arrangements is given by
$(2 - e^x)^{-1}$ and the
number of increasing preferential arrangements is $2^{n-1}$.

\begin{prop}
The quantity
$\left< \chi(\nabla(\e_n)), h_{1^n} \right>$ is a $q,t$
analog for the number of preferential arrangements (the maps
$f:[n] \rightarrow [k]$ where $1 \leq k \leq n$ which
are onto for some $k$).  Moreover, the quantity
$\left< \chi(\nabla(\e_n)), e_n \right>=\prod_{i=1}^{n-1}
(q^{n-i} + t^i)$ is a $q,t$ analog
of $2^{n-1}$.
\end{prop}

This proposition is a consequence of the
statement that appears in 
full generality just below.  For the moment we will provide
the following example:

\begin{example}
At $n=4$, there are $125=(4+1)^{4-1}$ parking functions, and 
$14=C_4$
weakly increasing parking functions represented by the following
list.  The first number is the number of parking 
functions such that
$(f(1),f(2),f(3),f(4))$ when sorted is the adjacent sequence.
The sum of these numbers is $125$.
\def\hs{{\hskip .1in}}
$$1 \times 1111 \hs
4 \times 1112 \hs
4 \times 1113 \hs
12 \times 1223 \hs
12 \times 1134 \hs
12 \times 1123 \hs
12 \times 1124 $$
$$4 \times 1222 \hs
4 \times 1114 \hs
6 \times 1133 \hs
6 \times 1122 \hs
12 \times 1224 \hs
12 \times 1233 \hs
24 \times 1234 $$
$75 = (2-e^x)^{-1}\coeff_{x^4} 4!$ of the parking
functions do not `skip' an integer, these are the 
preferential arrangements.  
Exactly $8 = 2^3$ of the preferential arrangements are
weakly increasing, those given by the following list:
$$1111\hs  1112 \hs 1122 \hs 1123\hs  1222 \hs
1223 \hs 1233 \hs1234$$

\end{example}

We remark that every preferential arrangement is also a
parking function.  This is a natural subset of the
parking functions which we will denote by ${Pref}_n$.
Just as in the case of the parking functions,
there is a natural $S_n$ action on this set formed by
permuting the values of the function (that is $(\sigma f)
(i) = f( \sigma_i )$) and hence ${Pref}_{n}$ forms
an $S_n$ module by defining a vector space with 
${ Pref}_n$ as the basis. 

Let $f$ be a preferential arrangement and define the content
of a preferential arrangement to be the composition 
$\alpha^{(f)}$ such that
the $i^{th}$ component is $|f^{-1}(i)|$.  We remark that two 
preferential arrangements are
in the same $S_n$-orbit if $\alpha^{(f)} = \alpha^{(g)}$.
The Frobenius series of
the $S_n$ module generated
by the preferential arrangements of content $\alpha$ is
given by the homogeneous symmetric function $h_\alpha$.

It follows that the preferential arrangement module may
be graded by a statistic on the content of the 
preferential arrangements.  If we choose our grading 
to be $q^{n(\inva{\alpha})}$ (this agrees with the 
`area' statistic on 
Dyck paths), then clearly the Frobenius series for the
module of preferential arrangements is given by
\begin{equation}
{\mathcal F}_{Pref_n}(q) = 
\sum_{\alpha \models n} q^{n(\inva{\alpha})} h_\alpha.
\end{equation}

In the commutative
case it is known that
\begin{equation}
\nabla( e_n )\coeff_{t=1} = \sum_{\mu \subseteq \delta_n} 
q^{{n\choose{2}}-|\mu|} e_{\la(\mu)},
\end{equation}
where $\delta_n = (n-1,n-2,\ldots,1,0)$ and $\la(\mu)$
is the sequence $(m_1(\mu), \ldots, m_{n-1}(\mu), 
n-\sum_{i=1}^{n-1} m_i(\mu))$ and $m_i(\mu)$ is the number
of parts of size $i$ in the partition $\mu$.  This
is related to the Frobenius series for the module of
parking functions by an application of the involution
$\omega$.
In exact analogy with the commutative case
we have the following proposition.

\begin{prop}
Set $t=1$ in the equation for the action of $\nabla$ on
$\e_n$, then
\begin{equation}
\nabla(\e_n) \coeff_{t=1} = 
\sum_{\alpha\models n} q^{n(\inva{\alpha})} \e_\alpha.
\end{equation}
\end{prop}

\begin{proof}With $t=1$ the action of $\nabla$
on $\s_{(1^n)}$ is given from equation (\ref{nabfor})
\begin{align}\nonumber
\nabla(\e_n) &= 
\sum_{\beta} 
\prod_{i \in D(\invb{\beta})} (1 + q^{i}) \s_\beta \\
&= 
\sum_{\beta} 
\sum_{\gamma \coarsereq \invb{\beta}} 
q^{n(\gamma)} \s_\beta \\\nonumber
&= 
\sum_{\gamma} 
\sum_{\invb{\beta} \finereq \gamma} 
q^{n(\gamma)} \s_\beta \\\nonumber
&= 
\sum_{\gamma} q^{n(\gamma)}
 \e_{\inva{\gamma}} = \sum_{\gamma} q^{n(\inva{\gamma})}
 \e_\gamma.
\end{align}

\end{proof}

This may be used to show that the quantity 
$\nabla( e_n )\coeff_{t=1} - \chi(\nabla(\e_n)) \coeff_{t=1}$
is $e$-positive (the coefficients are polynomials in
$q$ with non-negative integer coefficients when
the expression is expressed in the elementary basis.

We may use property (\ref{nabscal}) and (\ref{nabfor})
to calculate the 
inverse of this function as well.

\begin{prop}
\begin{equation}
\nabla^{-1}(\s_\alpha)=
(-1)^{\ell(\alpha)+1}\sum_{\invc{\alpha} \finereq \beta} 
q^{-n(\inva{\beta})}
t^{-n(\invc{\beta})}
\prod_{i \in D(\invc{\beta})\cap
D(\invc{\alpha})} (t^i + q^{n-i})
\s_\beta.
\end{equation}
\end{prop}

\begin{proof}
\begin{align}
\nabla^{-1}(\s_\alpha) &= 
\sum_{\beta\models n} (-1)^{n+\ell(\beta)}
\left< \nabla^{-1}(\s_\alpha), \s_{\invc{\beta}}\right>
\s_\beta \nonumber\\
&= \sum_{\beta\models n} (-1)^{n+\ell(\beta)} 
q^{-{n \choose 2}}
t^{-{n \choose 2}}
\left< \s_\alpha, \nabla(\s_{\invc{\beta}})\right>
\s_\beta\nonumber\\
&=\sum_{\beta\models n}\sum_{\gamma \finereq \beta}
(-1)^{\ell(\beta)+\ell(\invc{\beta})} 
q^{n(\invb{\beta})-{n \choose 2}}\\
&\hskip .2in t^{n(\beta)-{n \choose 2}}
\prod_{i \in D(\invc{\beta})\cap
D(\gamma)} (t^i + q^{n-i})
\left< \s_\alpha,  \s_\gamma \right>
\s_\beta\nonumber\\
&=\sum_{\invc{\alpha} \finereq \beta}
(-1)^{\ell(\alpha)+1} 
q^{-n(\inva{\beta})}
t^{-n(\invc{\beta})}
\prod_{i \in D(\invc{\beta})\cap
D(\invc{\alpha})} (t^i + q^{n-i})
\s_\beta.\nonumber
\end{align}
\end{proof}

\section{Appendix: Transition matrices between $\H_\alpha^q$ and
$\s_\beta$}

\begin{equation*}
\left[
\begin{array}{l|cc}
2&1&0\\
\noalign{\medskip}11&q&1
\end{array}
\right]
\end{equation*}

\begin{equation*}
\left[
\begin{array}{l|cccc}
3&1&0&0&0\\
\noalign{\medskip}12&q&1&0&0\\
\noalign{\medskip}21&q^2&0&1&0\\
\noalign{\medskip}111&{q}^{3}&q^2&{q}&1
\end{array}
\right]
\end{equation*}

\begin{equation*}
\left[\begin{array}{l|cccccccc} 
4&1&0&0&0&0&0&0&0\\
\noalign{\medskip}
13&q&1&0&0&0&0&0&0\\
\noalign{\medskip}22&{q}^{2}&0&1&0&0&0&0&0\\
\noalign{\medskip}112&{q}^{3}&{q}^{2}&q&1&0&0&0&0\\
\noalign{\medskip}31&{q}^{3}&0&0&0&1&0&0&0\\
\noalign{\medskip}121&{q}^{4}&{q}^{3}&0&0&q&1&0&0\\
\noalign{\medskip}211&{q}^{5}&0&{q}^{3}&0&{q}^{2}&0&1&0\\
\noalign{\medskip}
1111&{q}^{6}&{q}^{5}&{q}^{4}&{q}^{3}&{q}^{3}&{q}^{2}&q&1
\end{array}\right]
\end{equation*}

\section{Appendix: Transition matrices between $\Hw_\alpha^{qt}$
and $\s_\beta$}

\begin{equation*}
\left[\begin{array}{l|cc} 
2&1&q\\
\noalign{\medskip}
11&1&t
\end{array}
\right]
\end{equation*}

\begin{equation*}
\left[
\begin{array}{l|cccc}
3&1&{q}^{2}&q&{q}^{3}\\
\noalign{\medskip}
12&1&t&q&tq\\
\noalign{\medskip}
21&1&{q}^{2}&{t}^{2}&{t}^{2}{q}^{2}\\
\noalign{\medskip}
111&1&t&{t}^{2}&{t}^{3}
\end{array}\right]
\end{equation*}

\begin{equation*}
\left[
\begin{array}{l|cccccccc}
4&1&{q}^{3}&{q}^{2}&{q}^{5}&q&{q}^{4}&{q}^{3}&{q}^{6}\\
\noalign{\medskip}
13&1&t&{q}^{2}&t{q}^{2}&q&tq&{q}^{3}&t{q}^{3}\\
\noalign{\medskip}
22&1&{q}^{3}&{t}^{2}&{t}^{2}{q}^{3}&q&
{q}^{4}&{t}^{2}q&{t}^{2}{q}^{4}\\
\noalign{\medskip}
112&1&t&{t}^{2}&{t}^{3}&q&tq&{t}^{2}q&{t}^{3}q\\
\noalign{\medskip}
31&1&{q}^{3}&{q}^{2}&{q}^{5}&{t}^{3}&{t}^{3}{q}^{3}&
{t}^{3}{q}^{2}&{t}^{3}{q}^{5}\\
\noalign{\medskip}
121&1&t&{q}^{2}&t{q}^{2}&{t}^{3}&{t}^{4}&{t}^{3}{q}^{2}&
{t}^{4}{q}^{2}\\
\noalign{\medskip}
211&1&{q}^{3}&{t}^{2}&{t}^{2}{q}^{3}&{t}^{3}&
{t}^{3}{q}^{3}&{t}^{5}&{t}^{5}{q}^{3}\\
\noalign{\medskip}
1111&1&t&{t}^{2}&{t}^{3}&{t}^{3}&{t}^{4}&{t}^{5}&{t}^{6}
\end{array}\right]
\end{equation*}

\section{Appendix: Transition matrices between $\nabla(\s_\alpha)$
and $\s_\beta$}

\begin{equation*}
\left[\begin{array}{l|cc}
2&0&-qt\\
\noalign{\medskip}
11&1&t+q
\end{array}
\right]
\end{equation*}

\begin{equation*}
\left[
\begin{array}{l|cccc}
3&0&0&0&{q}^{3}{t}^{3}\\
\noalign{\medskip}
12&0&0&-q{t}^{2}&-\left (t+{q}^{2}\right )q{t}^{2}\\
\noalign{\medskip}
21&0&-{q}^{2}t&0&-\left ({t}^{2}+q\right ){q}^{2}t\\
\noalign{\medskip}
111&1&t+{q}^{2}&{t}^{2}+q&\left (t+{q}^{2}\right )
\left ({t}^{2}+q\right )
\end{array}
\right]
\end{equation*}

\vskip .5in
\noindent
{\bf Acknowledgement:}  Thank you to Andrew Rechnitzer and
Geanina Tudose for their valuable comments.  The authors
would also like to thank Fran\c cois Bergeron for 
discussions about the commutative operator $\nabla$.


\begin{thebibliography}{99}

\bibitem{A} E.~Abe, Hopf Algebras, 
New York, Cambridge University Press, 1980. 

\bibitem{BG} F.~Bergeron, A.~Garsia,
Science Fiction and Macdonald's Polynomials,
Algebraic Methods and $q$-Special Functions, 
R.~Floreanini, L.~Vinet (eds.), 
CRM Proceedings \& Lecture Notes,
American Mathematical Society, Volume 22, (1999), 1-52.

\bibitem{BGHT} F.~Bergeron, A.~Garsia,
M.~Haiman, G.~Tesler,
Identities and Positivity Conjectures for 
Some Remarkable Operators in the Theory of 
Symmetric Functions,  
Methods and Applications of Analysis, Volume 6, No. 3,  
(1999),  363-420. 

\bibitem{GH} A. Garsia, J. Haglund, A Proof of the 
$q,t$-Catalan 
Positivity Conjecture, preprint.
 
\bibitem{GKLLRT} I.~M.~Gelfand, D.~Krob, B.~Leclerc,
A.~Lascoux, V.~S.~Retakh, and J.-Y.~Thibon,
Noncommutative symmetric functions, 
Adv. Math. {\bf 112} (1995), 218-348.

\bibitem{H} M.~Haiman, Conjectures on the quotient ring 
by diagonal invariants, Journal of Algebraic 
Combinatorics \#3 (1994) 17-76.

\bibitem{Hi} F.~Hivert, Hecke Algebras, Difference Operators,
and Quasi-Symmetric Functions, Adv. Math. {\bf 155} (2000),
181-238.

\bibitem{HLT} F.~Hivert, A.~Lascoux, and J.-Y.~Thibon,
Noncommutative symmetric functions
and quasi-symmetric functions
with two and more parameters,
arXiv: {\tt math.CO/0106191}.

\bibitem{LLT} A.~Lascoux, B.~Leclerc, and J.-Y.~Thibon, 
Hecke algebras at roots of unity and crystal
bases of quantum affine algebras, Comm. 
Math. Phys. 181 (1997), 205-263.
 
\bibitem{M} I.~G.~Macdonald, Symmetric Functions and Hall
Polynomials, Oxford Mathematical Monographs, 
Oxford Univ. Press,
second edition, 1995.

\bibitem{MR} C.~Malvenuto and C.~Reutenauer, 
Duality between quasi-symmetric functions and the
Solomon descent algebra, {\it Journal of Algebra}, 
{\bf 177} (1995) 967--982.

\bibitem{S} M.~Sweedler, Hopf Algebras,
New York, W.A. Benjamin, 1969. 

\bibitem{St} R.~Stanley, Enumerative Combinatorics, Vol 1,
Wadsworth \& Brooks/Cole Advanced Books \& Software, 1986.

\bibitem{U} B.~C.~V.~Ung, NCSF, a Maple package for 
Noncommutative Symmetric Functions, Maple
Tech. News. 3, No. 3 (1996), 24-29.

\bibitem{UV} B.~C.~V.~Ung and S.~Veigneau, 
ACE une environnement 
en combinatoire alg\'ebrique,
in ``Proc. of the 7th Conf. Formal Power Series 
and Algebraic Combinatorics, 1995.''

\bibitem{Z} M.~Zabrocki, $q$-Analogs of Symmetric 
Function Operators, Discrete Mathematics (to appear).

\bibitem{Z2} M.~Zabrocki, Ribbon Operators and Hall-Littlewood Symmetric
Functions, Adv. in Math., {\bf 156} (2000), pp. 33-43.
       
\end{thebibliography}
\end{document}